\documentclass[12pt]{amsart}
 \usepackage{amsmath}
\usepackage{amssymb,latexsym}
\usepackage{amscd}
\usepackage{comment}
 \usepackage{xspace}
 \usepackage{verbatim}
\usepackage[english]{babel}

\begin{document}

\numberwithin{equation}{section}
\title[Positive solutions  of quasilinear elliptic equations]{Positive solutions  of quasilinear elliptic equations with subquadratic growth in the gradient}
\author{Moshe Marcus }
\address{Department of Mathematics, Technion\\
 Haifa 32000, ISRAEL}
 \email{marcusm@math.technion.ac.il}
\author{Phuoc-Tai Nguyen}
\address{Department of Mathematics, Technion\\
 Haifa 32000, ISRAEL}
 \email{nguyenphuoctai.hcmup@gmail.com}

\date{}

 \maketitle


\newcommand{\txt}[1]{\;\text{ #1 }\;}
\newcommand{\tbf}{\textbf}
\newcommand{\tit}{\textit}
\newcommand{\tsc}{\textsc}
\newcommand{\trm}{\textrm}
\newcommand{\mbf}{\mathbf}
\newcommand{\mrm}{\mathrm}
\newcommand{\bsym}{\boldsymbol}
\newcommand{\scs}{\scriptstyle}
\newcommand{\sss}{\scriptscriptstyle}
\newcommand{\txts}{\textstyle}
\newcommand{\dsps}{\displaystyle}
\newcommand{\fnz}{\footnotesize}
\newcommand{\scz}{\scriptsize}
\newcommand{\be}{\begin{equation}}
\newcommand{\bel}[1]{\begin{equation}\label{#1}}
\newcommand{\ee}{\end{equation}}
\newcommand{\eqnl}[2]{\begin{equation}\label{#1}{#2}\end{equation}}
\newcommand{\barr}{\begin{eqnarray}}
\newcommand{\earr}{\end{eqnarray}}
\newcommand{\bars}{\begin{eqnarray*}}
\newcommand{\ears}{\end{eqnarray*}}
\newcommand{\nnu}{\nonumber \\}
\newtheorem{subn}{\name}
\renewcommand{\thesubn}{}
\newcommand{\bsn}[1]{\def\name{#1}\begin{subn}}
\newcommand{\esn}{\end{subn}}
\newtheorem{sub}{\name}[section]
\newcommand{\dn}[1]{\def\name{#1}}   
\newcommand{\bs}{\begin{sub}}
\newcommand{\es}{\end{sub}}
\newcommand{\bsl}[1]{\begin{sub}\label{#1}}
\newcommand{\bth}[1]{\def\name{Theorem}
\begin{sub}\label{t:#1}}
\newcommand{\blemma}[1]{\def\name{Lemma}
\begin{sub}\label{l:#1}}
\newcommand{\bcor}[1]{\def\name{Corollary}
\begin{sub}\label{c:#1}}
\newcommand{\bdef}[1]{\def\name{Definition}
\begin{sub}\label{d:#1}}
\newcommand{\bprop}[1]{\def\name{Proposition}
\begin{sub}\label{p:#1}}
\newcommand{\R}{\eqref}
\newcommand{\rth}[1]{Theorem~\ref{t:#1}}
\newcommand{\rlemma}[1]{Lemma~\ref{l:#1}}
\newcommand{\rcor}[1]{Corollary~\ref{c:#1}}
\newcommand{\rdef}[1]{Definition~\ref{d:#1}}
\newcommand{\rprop}[1]{Proposition~\ref{p:#1}}
\newcommand{\BA}{\begin{array}}
\newcommand{\EA}{\end{array}}
\newcommand{\BAN}{\renewcommand{\arraystretch}{1.2}
\setlength{\arraycolsep}{2pt}\begin{array}}
\newcommand{\BAV}[2]{\renewcommand{\arraystretch}{#1}
\setlength{\arraycolsep}{#2}\begin{array}}
\newcommand{\BSA}{\begin{subarray}}
\newcommand{\ESA}{\end{subarray}}
\newcommand{\BAL}{\begin{aligned}}
\newcommand{\EAL}{\end{aligned}}
\newcommand{\BALG}{\begin{alignat}}
\newcommand{\EALG}{\end{alignat}}
\newcommand{\BALGN}{\begin{alignat*}}
\newcommand{\EALGN}{\end{alignat*}}
\newcommand{\note}[1]{\textit{#1.}\hspace{2mm}}
\newcommand{\Proof}{\note{Proof}}
\newcommand{\Remark}{\note{Remark}}
\newcommand{\modin}{$\,$\\[-4mm] \indent}
\newcommand{\forevery}{\quad \forall}
\newcommand{\set}[1]{\{#1\}}
\newcommand{\setdef}[2]{\{\,#1:\,#2\,\}}
\newcommand{\setm}[2]{\{\,#1\mid #2\,\}}
\newcommand{\mt}{\mapsto}
\newcommand{\lra}{\longrightarrow}
\newcommand{\lla}{\longleftarrow}
\newcommand{\llra}{\longleftrightarrow}
\newcommand{\Lra}{\Longrightarrow}
\newcommand{\Lla}{\Longleftarrow}
\newcommand{\Llra}{\Longleftrightarrow}
\newcommand{\warrow}{\rightharpoonup}
\newcommand{
\paran}[1]{\left (#1 \right )}
\newcommand{\sqbr}[1]{\left [#1 \right ]}
\newcommand{\curlybr}[1]{\left \{#1 \right \}}
\newcommand{\abs}[1]{\left |#1\right |}
\newcommand{\norm}[1]{\left \|#1\right \|}
\newcommand{
\paranb}[1]{\big (#1 \big )}
\newcommand{\lsqbrb}[1]{\big [#1 \big ]}
\newcommand{\lcurlybrb}[1]{\big \{#1 \big \}}
\newcommand{\absb}[1]{\big |#1\big |}
\newcommand{\normb}[1]{\big \|#1\big \|}
\newcommand{
\paranB}[1]{\Big (#1 \Big )}
\newcommand{\absB}[1]{\Big |#1\Big |}
\newcommand{\normB}[1]{\Big \|#1\Big \|}
\newcommand{\produal}[1]{\langle #1 \rangle}

\newcommand{\thkl}{\rule[-.5mm]{.3mm}{3mm}}
\newcommand{\thknorm}[1]{\thkl #1 \thkl\,}
\newcommand{\trinorm}[1]{|\!|\!| #1 |\!|\!|\,}
\newcommand{\bang}[1]{\langle #1 \rangle}
\def\angb<#1>{\langle #1 \rangle}
\newcommand{\vstrut}[1]{\rule{0mm}{#1}}
\newcommand{\rec}[1]{\frac{1}{#1}}
\newcommand{\opname}[1]{\mbox{\rm #1}\,}
\newcommand{\supp}{\opname{supp}}
\newcommand{\dist}{\opname{dist}}
\newcommand{\myfrac}[2]{{\displaystyle \frac{#1}{#2} }}
\newcommand{\myint}[2]{{\displaystyle \int_{#1}^{#2}}}
\newcommand{\mysum}[2]{{\displaystyle \sum_{#1}^{#2}}}
\newcommand {\dint}{{\displaystyle \myint\!\!\myint}}
\newcommand{\q}{\quad}
\newcommand{\qq}{\qquad}
\newcommand{\hsp}[1]{\hspace{#1mm}}
\newcommand{\vsp}[1]{\vspace{#1mm}}
\newcommand{\ity}{\infty}
\newcommand{\prt}{\partial}
\newcommand{\sms}{\setminus}
\newcommand{\ems}{\emptyset}
\newcommand{\ti}{\times}
\newcommand{\pr}{^\prime}
\newcommand{\ppr}{^{\prime\prime}}
\newcommand{\tl}{\tilde}
\newcommand{\sbs}{\subset}
\newcommand{\sbeq}{\subseteq}
\newcommand{\nind}{\noindent}
\newcommand{\ind}{\indent}
\newcommand{\ovl}{\overline}
\newcommand{\unl}{\underline}
\newcommand{\nin}{\not\in}
\newcommand{\pfrac}[2]{\genfrac{(}{)}{}{}{#1}{#2}}

\def\ga{\alpha}     \def\gb{\beta}       \def\gg{\gamma}
\def\gc{\chi}       \def\gd{\delta}      \def\ge{\epsilon}
\def\gth{\theta}                         \def\vge{\varepsilon}
\def\gf{\phi}       \def\vgf{\varphi}    \def\gh{\eta}
\def\gi{\iota}      \def\gk{\kappa}      \def\gl{\lambda}
\def\gm{\mu}        \def\gn{\nu}         \def\gp{\pi}
\def\vgp{\varpi}    \def\gr{\rho}        \def\vgr{\varrho}
\def\gs{\sigma}     \def\vgs{\varsigma}  \def\gt{\tau}
\def\gu{\upsilon}   \def\gv{\vartheta}   \def\gw{\omega}
\def\gx{\xi}        \def\gy{\psi}        \def\gz{\zeta}
\def\Gg{\Gamma}     \def\Gd{\Delta}      \def\Gf{\Phi}
\def\Gth{\Theta}
\def\Gl{\Lambda}    \def\Gs{\Sigma}      \def\Gp{\Pi}
\def\Gw{\Omega}     \def\Gx{\Xi}         \def\Gy{\Psi}

\def\CS{{\mathcal S}}   \def\CM{{\mathcal M}}   \def\CN{{\mathcal N}}
\def\CR{{\mathcal R}}   \def\CO{{\mathcal O}}   \def\CP{{\mathcal P}}
\def\CA{{\mathcal A}}   \def\CB{{\mathcal B}}   \def\CC{{\mathcal C}}
\def\CD{{\mathcal D}}   \def\CE{{\mathcal E}}   \def\CF{{\mathcal F}}
\def\CG{{\mathcal G}}   \def\CH{{\mathcal H}}   \def\CI{{\mathcal I}}
\def\CJ{{\mathcal J}}   \def\CK{{\mathcal K}}   \def\CL{{\mathcal L}}
\def\CT{{\mathcal T}}   \def\CU{{\mathcal U}}   \def\CV{{\mathcal V}}
\def\CZ{{\mathcal Z}}   \def\CX{{\mathcal X}}   \def\CY{{\mathcal Y}}
\def\CW{{\mathcal W}} \def\CQ{{\mathcal Q}}
\def\BBA {\mathbb A}   \def\BBb {\mathbb B}    \def\BBC {\mathbb C}
\def\BBD {\mathbb D}   \def\BBE {\mathbb E}    \def\BBF {\mathbb F}
\def\BBG {\mathbb G}   \def\BBH {\mathbb H}    \def\BBI {\mathbb I}
\def\BBJ {\mathbb J}   \def\BBK {\mathbb K}    \def\BBL {\mathbb L}
\def\BBM {\mathbb M}   \def\BBN {\mathbb N}    \def\BBO {\mathbb O}
\def\BBP {\mathbb P}   \def\BBR {\mathbb R}    \def\BBS {\mathbb S}
\def\BBT {\mathbb T}   \def\BBU {\mathbb U}    \def\BBV {\mathbb V}
\def\BBW {\mathbb W}   \def\BBX {\mathbb X}    \def\BBY {\mathbb Y}
\def\BBZ {\mathbb Z}

\def\GTA {\mathfrak A}   \def\GTB {\mathfrak B}    \def\GTC {\mathfrak C}
\def\GTD {\mathfrak D}   \def\GTE {\mathfrak E}    \def\GTF {\mathfrak F}
\def\GTG {\mathfrak G}   \def\GTH {\mathfrak H}    \def\GTI {\mathfrak I}
\def\GTJ {\mathfrak J}   \def\GTK {\mathfrak K}    \def\GTL {\mathfrak L}
\def\GTM {\mathfrak M}   \def\GTN {\mathfrak N}    \def\GTO {\mathfrak O}
\def\GTP {\mathfrak P}   \def\GTR {\mathfrak R}    \def\GTS {\mathfrak S}
\def\GTT {\mathfrak T}   \def\GTU {\mathfrak U}    \def\GTV {\mathfrak V}
\def\GTW {\mathfrak W}   \def\GTX {\mathfrak X}    \def\GTY {\mathfrak Y}
\def\GTZ {\mathfrak Z}   \def\GTQ {\mathfrak Q}

\def\sign{\mathrm{sign\,}}
\def\bdw{\prt\Gw\xspace}
\def\nabu{|\nabla u|}


\begin{abstract}
We study positive solutions of equation (E) $-\Gd u + u^p|\nabla u|^q= 0$ ($0<p$, $0\leq q\leq 2$, $p+q>1$) and other related equations  in a smooth bounded domain $\Gw\sbs\BBR^N$.  We show that if $N(p+q-1)<p+1$ then, for every positive, finite Borel measure $\mu$ on $\bdw$,  there exists a solution of (E) such that $u=\mu$ on $\prt \Gw$. Furthermore, if $N(p+q-1)\geq p+1$ then an isolated point singularity on $\bdw$ is removable. In particular there is no solution with boundary data $\gd_y$ (=Dirac measure at a point $y\in \bdw$). Finally we obtain a classification of positive solutions with an isolated boundary singularity.
\end{abstract}

\noindent {\bf Keywords:}  quasilinear equation;  boundary singularities; Radon measures; weak singularities; strong singularities; boundary trace, removability.

\medskip

\section{Introduction}
In this paper, we are concerned with the boundary data measure problem associated to the equation
\bel{A0} -\Gd u + H(x,u,\nabla u) = 0 \ee
in $\Gw$ where $\Gw$ is a domain in $\BBR^N$ and $H$ is a Caratheodory function defined in $\Gw \ti \BBR \ti \BBR^N$.

When $H$ depends only on $u$, much works on the boundary value problem for equation \eqref{A0}, especially for the following typical equation
\bel{Ap} -\Gd u + u^p = 0  \ee
with $p>1$, have been studied by Le Gall \cite{Lg2} , Gmira and V\'eron \cite{GV}, Marcus and V\'eron \cite{MV1}, \cite{MV2}, \cite{V1}, \cite{V2}. It was shown that equation \eqref{Ap} admits a {\it critical value}
\bel{pc} p_c=\frac{N+1}{N-1}. \ee
For any $1<p<p_c$, if $\gm$ is a bounded Radon measure on $\prt \Gw$, then there exists a unique solution of \eqref{Ap} with boundary data $\gm$. Moreover isolated boundary singularities of solutions of \eqref{Ap} can be completely described. More precisely, when $1<p<p_c$, if $u \in C^2(\Gw) \cap C(\ovl \Gw \sms \{0\})$ is a nonnegative solution of  \eqref{Ap} vanishing on $\prt \Gw \sms \{0\}$  then either $u$ behaves like $kP^\Gw(.,0)$ near $0$ with some $k \geq 0$ and $P^\Gw(.,0)$ being the Poisson kernel in $\Gw$, or $u(x) \approx \rho(x) \abs{x}^{-\frac{q+1}{q-1}}$ as $x \to 0$ where $\rho$ is the distance function to $\prt \Gw$. When $p \geq p_c$, the isolated singularities are {\it  removable}, namely if $u \in C^2(\Gw) \cap C(\ovl \Gw \sms \{0\})$ is a nonnegative solution of  \eqref{Ap} vanishing on $\prt \Gw \sms \{0\}$ then $u \equiv 0$. More general results can be found in \cite{MV3}, \cite{MV4}.

The case $H$ depends only on $\nabla u$ has been recently investigated by P.T. Nguyen and L. V\'eron \cite{NV}. Many results have been extended to quasilinear equations of the form
	\bel{Aq} -\Gd u + g(|\nabla u|) = 0 \ee
in $\Gw$. Under suitable conditions on $g$, if $\gm$ is a bounded Radon measure on $\prt \Gw$, they proved existence of a positive solution of \eqref{Aq} with boundary data $\gm$. In the power case, namely $g(|\nabla u|)=|\nabla u|^q$ with $1 \leq q \leq 2$, they showed that the critical value for \eqref{Aq} is
\bel{qc} q_c=\frac{N+1}{N} \ee
and analogous phenomena occur for isolated boundary singularities. Notice that when $q >2$, by \cite{Li1}  if $u \in C^2(\Gw)$ is a positive solution of \eqref{Aq} then $u$ is bounded in $\Gw$, therefore there is no singularity on  the boundary.

Motivated by the above papers, we study boundary singularities of positive solutions of \eqref{A0} in the case that {\it $H$ depends on both $u$ and $\nabla u$}. We are interested in the case of subquadratic growth in the gradient and concentrate in particular on two model cases
\bel{multi} H(x,u,\gx) =   u^p|\gx|^q \forevery (x,u,\gx) \in \Gw \ti \BBR_+ \ti \BBR^N \ee
where $p>0$, $0 \leq q \leq 2$ and
\bel{add}  H(x,u,\gx)  = u^p + |\gx|^q  \forevery (x,u,\gx) \in \Gw \ti \BBR_+ \ti \BBR^N \ee
where $p \geq 1$, $1 \leq q \leq 2$. Concerning the above types of nonlinearity, there have been many works on {\it large solutions}, namely solutions that blow up on the boundary. When $H$ satisfies \eqref{multi}, there exists no large solution to equation \eqref{A0}. When $H$ satisfies \eqref{add} there exists a large solution to \eqref{A0}; moreover large solution is unique if $1<p<q \leq 2$ (see \cite{AGMQ}, \cite{BaGi}). To our knowledge, up to now, no study dealing with the boundary value problem with measure data for these types of nonlinearity has been published. We list below results concerning existence of solution with boundary data as Radon measure, classification of isolated boundary singularities in subcritical case and removability in critical and supercritical case.

In what follows, unless otherwise stated, $\Gw$ is a bounded domain of class $C^2$ with $\prt \Gw$ containing the origin $0$, $S^{N-1}$ the unit sphere, $S_+^{N-1}=S^{N-1} \cap \BBR_+^{N}$ the upper hemisphere and $(r,\gs) \in \BBR_+ \ti S^{N-1}$ the spherical coordinates in $\BBR^N $.  To state our main results, it is convenient to introduce the definition of solutions.
\bdef{sol} i) A function $u$ is called a solution of \eqref{A0} if $u \in L_{loc}^1(\Gw)$, $H(x,u,\nabla u) \in L_{loc}^1(\Gw)$ and $u$ satisfies \eqref{A0} in the sense of distribution, i.e.
	$$\myint{\Gw}{}\left(-u\Gd\gz+H(x,u,\nabla u) \gz\right)dx= 0 $$
for every $\gz \in C_c^\infty(\Gw)$. \medskip

\noindent ii) Let $\gm$ is a positive Borel measure on $\prt\Gw$. A function $u$ is called a solution of
\bel{A1}\left\{\BA {l} -\Gd u +  H(x,u,\nabla u) = 0 \qq \text{in } \Gw\\[1mm]
\phantom{-\Gd  +g(\abs{\nabla u}),,,,}
u=\gm \qq \text{on } \prt\Gw
\EA\right.\ee
if $u \in L^1(\Gw)$, $H(x,u,\nabla u) \in L^1_\rho(\Gw)$ where $\rho(x):=\dist (x,\prt\Gw)$ and $u$ satisfies
\bel{A2}
\myint{\Gw}{}\left(-u\Gd\gz+H(x,u,\nabla u)\gz\right)dx=-\myint{\prt\Gw}{}\myfrac{\prt \gz}{\prt \bf n}d\gm
\ee
for all  $\gz \in C_0^2(\ovl \Gw)$, where $\bf n$ denotes the normal outward unit vector to $\prt \Gw$.
\es

\bdef{subcriticality} A nonlinearity $H$ is called {\bf subcritical} if the problem \eqref{A1} admits a solution for every positive bounded measure $\gm$ on $\prt \Gw$. Otherwise, $H$ is called {\bf supercritical}. \es
%
%
%
Set
\bel{mpq} m_{p,q}=\max\left\{p,\frac{q}{2-q}\right\}. \ee
Following is the main existence result in the subcritical case. \medskip

\noindent{\bf Theorem A.} {\it Assume either $H$ satisfies \eqref{multi} with $0<N(p+q-1)<p+1$ or $H$ satisfies \eqref{add} with $m_{p,q}<p_c$. Then $H$ is subcritical. Moreover, let $\{\gm_n\}$ be a sequence of positive bounded measures on $\prt\Gw$ which converges to a positive bounded $\gm$ in the weak sense of measures and $\{u_{\gm_n}\}$ be a sequence of corresponding solutions of \eqref {A1} with $\gm=\gm_n$. Then there exists a subsequence such that $\{ u_{\gm_{n_k}}\}$ converges to a solution $ u_\gm$ of \eqref {A1} in $L^1(\Gw)$ and $\{H(x,u_{\gm_{n_k}},\nabla u_{\gm_{n_k}})\}$ converges to $H(x,u,\nabla u)$ in $L^1_\rho(\Gw)$.} \medskip

\noindent {\bf Remark.} The method used is classical, using the estimates in weak $L^p$ space and compactness of approximating solutions. Due to this approach, the results stated in Theorem A can be extended to the following cases:
\bel{multia} 0 \leq H(x,u,\gx) \leq a_1(x)u^p|\gx|^q  \forevery (x,u,\gx) \in \Gw \ti \BBR_+ \ti \BBR^N \ee
where $p>0$, $q \geq 0$,  $0<N(p+q-1)<p+1$ , $a_1 \in L^\infty(\Gw)$ and $a_1>c>0$;
\bel{adda} 0 \leq H(x,u,\nabla u)  \leq a_2(x)f(u) + a_3(x)g(|\gx|)  \forevery (x,u,\gx) \in \Gw \ti \BBR_+ \ti \BBR^N \ee
where $a_i \in L^\infty(\Gw)$, $a_i>c>0$ ($i=3,4$), $f$ and $g$ are positive, nondecreasing, continuous functions in $\BBR_+$, satisfying $f(0)=g(0)=0$ and
$$ \myint{1}{\infty}t^{-\frac{2N}{N-1}}f(t)dt<\infty, \q \myint{1}{\infty}t^{-\frac{2N+1}{N}}g(t)dt<\infty. $$

The uniqueness of the problem remains open. However, if  $\gm$ is concentrated at a point on the boundary and the functions $a_i$ ($i=1,2,3$) are positive constants, we prove that the solution of \eqref{A1} is unique. \medskip

\noindent{\bf Theorem B.} {\it Assume either $H$ satisfies \eqref{multi} with $0<N(p+q-1)<p+1$ or $H$ satisfies \eqref{add} with $m_{p,q}<p_c$. Then for any $k>0$, there exists a unique positive solution to \eqref{A1} with $\gm=k\gd_0$,  denoted by $u^\Gw_{k,0}$, where $\gd_0$ is the Dirac mass concentrated at the origin $0$. Moreover,
\bel{weaksing} u^{\Gw}_{k,0}(x)=kP^\Gw(x,0)(1+o(1)) \q \text{as } x \to 0. \ee
and there exists $d_k>0$ such that
\bel{com} d_k P^\Gw(x,0)<u(x) < k P^\Gw(x,0) \forevery x \in \Gw, \ee
}

The solutions $u_{k,0}^\Gw$ are called {\bf weakly singular solutions}. It follows from \eqref{weaksing} that the sequence $\{u_{k,0}^\Gw\}$ is increasing.
Hence, it is interesting to study the limit of this sequence. In order to state the result involving the limit, we define the class of {\bf strongly singular solutions} (see the definition of the boundary trace $tr_{\prt \Gw}$ in section 3.2)
\bel{strongset} \CU^\Gw_0:=\{u \in C^2(\Gw) \text{ positive solution of \eqref{A0} with } tr_{\prt \Gw}(u)=(\{0\},0)  \}. \ee

\noindent{\bf Theorem C.} {\it Under the assumptions of theorem B, the function $u^\Gw_{\infty,0}:=\lim_{k \to \infty}u^\Gw_{k,0}$ is the minimal element of $\CU^\Gw_0$.} \medskip

The asymptotic behavior of $u^\Gw_{\infty,0}$ can be obtained  due to the study solution of
	\bel{P3} \left\{ \BA{lll}
		- \Gd u + H(x,u,\nabla u) = 0 \qq &\text{ in } \BBR^N_+ \\
		\phantom{- \Gd  + \abs{\nabla u}^qu^p,,,,}
		u = 0 &\text{ on } \prt \BBR^N_+\sms\{0\}
	\EA \right. \ee
under the separate form $u(x)=r^{-\gb}\gw(\gs)$ where $\gb>0$, $r=|x|$ and $\gs=\frac{x}{|x|} \in S^{N-1}_+$. Denote by $\nabla' $ and  $\Gd'$ the covariant derivative on $S^{N-1}$ identified with the tangential derivative and the Laplace-Beltrami operator on $S^{N-1}$ respectively.

When $H$ satisfies \eqref{multi}, by plugging $u=r^{-\gb}\gw(\gs)$ into \eqref{P3} we deduce that
\bel{beta1}\gb=\gb_1:=\frac{2-q}{p+q-1} \ee
 and $\gw$ satisfies
\bel{PH1} -\Gd' \gw +  F_1(\gw,\nabla' \gw)=0 \text{ in } S_+^{N-1}, \q \gw = 0 \text{ on } \prt S_+^{N-1} \ee
where $ F_1(s,\gx):=  s^p(\gb_1^2\, s^2+\abs{\gx}^2)^{\frac{q}{2}}-\gb_1(\gb_1+2-N)s $ with $s \in \BBR_+$ and $\gx \in \BBR^N$.

When $H$ satisfies \eqref{add}  we deduce that
\bel{beta2}\gb=\gb_2:=\frac{2}{m_{p,q}-1} \ee
where $m_{p,q}$ is defined in \eqref{mpq}. Moreover if $p=\frac{q}{2-q}$ then $\gw$ satisfies
\bel{PH2} -\Gd' \gw +  F_2(\gw,\nabla' \gw)=0 \text{ in } S_+^{N-1}, \q \gw = 0 \text{ on } \prt S_+^{N-1} \ee
where  $ F_2(s,\gx):= s^p + \ga_3(\gb_{2}^2\, s^2+\abs{\gx}^2)^{\frac{q}{2}}-\gb_{2}(\gb_{2}+2-N)s $ with $s \in \BBR_+$ and $\gx \in \BBR^N$. When $p>\frac{q}{2-q}$, we consider
\bel{PH3} -\Gd' \gw +  F_3(\gw, \nabla' \gw)=0 \text{ in } S_+^{N-1}, \q \gw = 0 \text{ on } \prt S_+^{N-1} \ee
where  $ F_3(s,\gx):= s^p-\gb_{2}(\gb_{2}+2-N)s $ with $s \in \BBR_+$. When $p<\frac{q}{2-q}$, we consider
\bel{PH4} -\Gd' \gw +  F_4(\gw,\nabla' \gw)=0 \text{ in } S_+^{N-1}, \q \gw = 0 \text{ on } \prt S_+^{N-1} \ee
where  $ F_4(s,\gx):=\ga_3(\gb_{2}^2\, s^2+\abs{\gx}^2)^{\frac{q}{2}}-\gb_{2}(\gb_{2}+2-N)s $ with $s \in \BBR_+$ and $\gx \in \BBR^N$.

Denote by $\CE_i$ $(i=\ovl {1,4}$) the set of positive solutions in $C^2(S^{N-1}_+)$ of
\bel{PHi} -\Gd' \gw +  F_i(\gw,\nabla' \gw)=0 \text{ in } S_+^{N-1}, \q \gw = 0 \text{ on } \prt S_+^{N-1}. \ee
\medskip

\noindent{\bf Theorem D.} {\it i) If $H$ satisfies \eqref{multi} with $0<N(p+q-1)<p+1$ then $\CE_1 \neq \emptyset$. Moreover, if $p \geq 1$ then there exists a unique solution $\gw^*_1$ of \eqref{PH1}, namely $\CE_1=\{\gw^*_1\}$. In addition, $\CU^\Gw_0=\{u^\Gw_{\infty,0}\}$ and
\bel{E1}
\lim_{\tiny\BA{c}\Gw \ni x\to 0\\
\frac{x}{|x|}=\gs\in S^{N-1}_+
\EA}|x|^{\gb_1}u^\Gw_{\infty,0}(x)=\gw^*_1(\gs)
\ee
 locally uniformly on $S^{N-1}_+$. \smallskip

\noindent ii) If $H$ satisfies \eqref{add} with $m_{p,q}<p_c$ then there $\CE_i = \{\gw_i^*\}$ where $i=2$ if $p=\frac{q}{2-q}$, $i=3$ if $p>\frac{q}{2-q}$, $i=4$ if $p<\frac{q}{2-q}$. In addition, $\CU^\Gw_0=\{u^\Gw_{\infty,0}\}$ and
\bel{E1'}
\lim_{\tiny\BA{c}\Gw \ni x\to 0\\
\frac{x}{|x|}=\gs\in S^{N-1}_+
\EA}|x|^{\gb_2}u^\Gw_{\infty,0}(x)=\gw_i^*(\gs)
 \ee
 locally uniformly on $S^{N-1}_+$ .} \medskip

\noindent {\bf Remark.} Notice that when $H$ satisfies \eqref{multi} or $H$ satisfies \eqref{add} with $p=\frac{q}{2-q}$, the equation \eqref{A0} is unvariant under an appropriate similarity transformation. However, it is not that case when $H$ satisfies \eqref{add} with $p \neq \frac{q}{2-q}$; in this situation, there is a competition between $u^p$ and $|\nabla u|^q$. The theorem D shows that when $p>\frac{q}{2-q}$, the term $u^p$ plays a dominant role and hence the solution $u^\Gw_{\infty,0}$ behaves like $u^\Gw_{p,\infty,0}$ near $0$ where  $u^\Gw_{p,\infty,0}$ is the solution of
\bel{up-infty} -\Gd u +   u^p = 0 \text{ in } \Gw, \q tr_{\prt \Gw}(u)=(\{0\},0). \ee
Otherwise, when $p<\frac{q}{2-q}$, $|\nabla u|^q$ is the dominant term and therefore $u^\Gw_{\infty,0}$ behaves like $u^\Gw_{q,\infty,0}$ near $0$ where  $u^\Gw_{q,\infty,0}$ is the solution of
\bel{uq-infty} -\Gd u + |\nabla u|^q= 0 \text{ in } \Gw, \q tr_{\prt \Gw}(u)=(\{0\},0). \ee

As a consequence, we provide a full characterization of isolated singularities at the origin $0$. \medskip

\noindent{\bf Theorem E} {\it Assume either $H$ satisfies \eqref{multi} with $0<N(p+q-1)<p+1$ and $p \geq 1$ or $H$ satisfies \eqref{add} with $m_{p,q}<p_c$. Let $u \in C(\ovl \Gw \sms \{0\}) \cap C^2(\Gw)$ be a nonnegative solution of \eqref{A0} vanishing on $\prt \Gw \sms \{0\}$. Then
\begin{itemize}
\item either $u \equiv 0$,
\item or there exists $k>0$ such that $u(x)=u^{\Gw}_{k,0}=kP^\Gw(x,0)(1+o(1))$ as $x \to 0$,
\item or $u(x)=u^{\Gw}_{\infty,0}$, the unique element of $\CU^\Gw_0$, and the asymptotic behavior of $u$ near $0$ is given either in \eqref{E1} or in \eqref{E1'} according to the assumptions on $H$.
\end{itemize}
 } \medskip

On the contrary, we show that  isolated boundary singularities are {\bf removable} in the critical and supercritical case. More precisely, \medskip

\noindent{\bf Theorem F} {\it Assume either $H$ satisfies \eqref{multi} with $N(p+q-1) \geq p+1$
or H satisfies \eqref{add} with $m_{p,q} \geq p_c$. If $u\in C(\overline\Gw\setminus\{0\})\cap C^2(\Gw)$ is a nonnegative solution of \eqref{A0} vanishing on $\prt\Gw\setminus\{0\}$ then $u\equiv 0$.} \medskip

When $H$ satisfies \eqref{multi}, the proof of Theorem F is divided into three cases. The case $N(p+q-1)>p+1$ is treated due to a priori estimate for solutions with isolated singularity at $0$. The critical case is more delicate: we first prove removability result for $\Gw=\BBR^N_+$ and then by using regularity results up to boundary (see \cite{Lib}) we get the assertion when $\Gw$ is bounded. Finally, when $q=2$, thanks to a change of unknown, we deduce that $u \equiv 0$. When $H$ satisfies \eqref{add}, the removabilty result for \eqref{A0} is derived from the one for \eqref{Ap} and \eqref{Aq}.

The paper is organized as follows. In section 2, we establish some estimates on positive solution of \eqref{A0} and its gradient, and recall some estimates concerning weak $L^p$ space which play a key role in proving the existence of solutions with bounded boundary measure data in the subcritical case. Section 3 is devoted to the proof of Theorem A and to investigate the notion of {\it boundary trace}. In section 4, we provide a complete description of isolated singularities (Theorem B, Theorem C, Theorem D and Theorem E). Finally, in section 5, we give proof of removability result (Theorem F).

Throughout the present paper, we denote by $c$, $C$, $c_1$, $c_2$,...positive constants which may vary from line to line. If necessary the dependence of these constants will be made precise.

\section{Preliminaries}
The following comparison principle can be found in \cite[Theorem 9.2]{GT}.
\bprop{comparison} Assume $H: \Gw \ti \BBR_+ \ti \BBR^N \to \BBR_+$ is nondecreasing with respect to $u$ for any $(x,\gx) \in \Gw \ti \BBR^N$, continuously differentiable with respect to $\gx$ and $H(x,0,0)=0$. Let $u_1$, $ u_2 \in C^2(\Gw) \cap C(\ovl \Gw)$ be two nonnegative solution of \eqref{A0}. If
$$ -\Gd u_1 + H(x,u_1,\nabla u_1) \leq -\Gd u_2 + H(x,u_2,\nabla u_2) \quad \text{in } \Gw $$
and $u_1 \leq u_2$ on $\prt \Gw$. Then $u_1 \leq u_2$ in $\Gw$.
\es

Next, for $\gd>0$, we set
	$$ \Gw_\gd=\{x \in \Gw: \rho(x) < \gd\},\q D_\gd=\{x \in \Gw: \rho(x) > \gd\}, $$
	$$ \Gs_\gd=\prt D_\gd=\{x \in \Gw: \rho(x) = \gd\},\q \Gs=\prt \Gw. $$
Since $\Gw$ is of class $C^2$, there exists $\gd_0>0$ such that \medskip

\noindent i) For every $x \in \ovl \Gw_{\gd_0}$, there exists a unique point $\gs(x) \in \prt \Gw$ such that $x=\gs(x)-\rho(x) {\bf n}_{\gs(x)}$ where ${\bf n}_{\gs(x)}$ is the outward unit normal vector to $\prt \Gw$ at $\gs(x)$. \medskip

\noindent ii) The mappings $x \mapsto \rho(x)$ and $x \mapsto \gs(x)$ belong to $C^2(\ovl \Gw_{\gd_0})$ and $C^1(\ovl \Gw_{\gd_0})$ respectively. Moreover, $\lim_{x \to \gs(x)}\nabla \rho(x)=-{\bf n}_{\gs(x)}$ and $\abs{\nabla \rho}=1$ in $\Gw_{\gd_0}$. \medskip

As a consequence of \rprop{comparison}, we deduce  a priori estimate
\bprop{est-multi} Assume $H$ satisfies \eqref{multi} with $p \geq 0$, $0\leq q <2$, $p+q>1$. Let $u \in C^2(\Gw)$ be a positive solution of equation \eqref{A0}. Then
\bel{est1} u(x) \leq \Gl_1\,\rho(x)^{-\gb_1} +\max\{u(x): x \in \ovl D_{\gd_0}\}, \quad \forall x \in \Gw, \ee
\bel{est2} u(x) \leq  \Gl_1\rho(x)^{-\gb_1} + \Gl_1'\norm{u}_{L^1(D_{\frac{\gd_0}{2}})} \quad \forall x \in \Gw,\ee
\bel{grad} \abs{\nabla u(x)} \leq \tl \Gl_1\,\rho(x)^{-\gb_1-1} \forevery x \in \Gw\ee
where $\gb_1$ is defined in \eqref{beta1}, $\Gl_1'=\Gl_1'(N,\gd_0)$,
$\tl \Gl_1=\tl \Gl_1(N,p,q,\Gw,\norm{u}_{L^1(D_{\frac{\gd_0}{2}})})$  and
\bel{Lambda1}  \Gl_1=\left(\myfrac{\gb_1+2}{ \gb_1^{q-1}}\right)^{\frac{1}{p+q-1}}. \ee
\es
\Proof {\it Proof of \eqref{est1}.} Put $ M_{\gd_0}=\max\{u(x): x \in \ovl D_{\gd_0}\}$. For each $\gd \in (0,\gd_0)$, we set $w_{\gd}(x)=\Gl_1(\rho(x)-\gd)^{-\gb_1}+M_{\gd_0}$ for $x \in D_\gd$. We can choose $\gd_0< \norm{\Gd \rho}^{-1}_{L^\infty(\Gw)}$. By a computation, we obtain in $\Gw_{\gd_0}\sms \ovl \Gw_{\gd}$,
$$  -\Gd w_{\gd} +  w_{\gd}^p\abs{\nabla w_{\gd}}^q > 0. $$
Since $w_\gd \geq u$ on $\Gs_{\gd} \cup \Gs_{\gd_0}$, by comparison principle \rprop{comparison}, $u \leq w_\gd$ in $\Gw_{\gd_0} \sms \ovl \Gw_{\gd}$. Letting $\gd \to 0$ leads to the conclusion. \medskip

\noindent {\it Proof of \eqref{est2}.} The estimate \eqref{est2} follows from \eqref{est1} and \cite[Theorem 1.1]{Tr1}.  \medskip

\noindent{\it Proof of \eqref{grad}.} Fix $x_0 \in \Gw_{\frac{3\gd_0}{4}}$ and set
$$d_0=\frac{1}{3}\rho(x_0), \quad M_{0}=\max\{u(x): x \in B_{2d_0}(x_0)\}, $$
$$ \quad u_0(y)=\myfrac{u(x)}{M_0},\quad  y=\frac{1}{d_0}x \in B_2(y_0),\quad  y_0=\frac{1}{d_0}x_0. $$
Then $\max\{u_0(y): y \in B_2(y_0)\}=1$ and $u_0$ satisfies
	$$ -\Gd u_0 +   M_0^{p+q-1}d_0^{2-q}u_0^p|\nabla u_0|^q = 0 $$
in $B_2(y_0)$. It follows from \cite{La} that there exists a positive constant $c=c(N,p,q,\gd_0,\norm{u}_{L^1(D_{\frac{\gd_0}{2}})})$ such that $\max_{B_1(y_0)}|\nabla u_0| \leq c$. Consequently,
	$$ \max_{B_1(y_0)}|\nabla u| \leq \myfrac{c}{d_0}\max_{B_2(y_0)}u. $$
Therefore, we deduce \eqref{grad}. \qed \medskip

By an analogous argument, we obtain

\bprop{est-add} Assume $H$ satisfies \eqref{add} with $p > 1$, $1< q <2$. Let $u \in C^2(\Gw)$ be a positive solution of equation \eqref{A0}. Then
\bel{est2-add} u(x) \leq  \Gl_2\rho(x)^{-\gb_2} + \Gl_2'\norm{u}_{L^1(D_{\frac{\gd_0}{2}})} \ee
\bel{grad-add} \abs{\nabla u(x)} \leq \tl \Gl_2\,\rho(x)^{-\gb_2-1} \forevery x \in \Gw\ee
where $\gb_2$ is defined in \eqref{beta2}, $\Gl_2=\Gl_2(p, q)$, $\Gl_2'=\Gl_2'(N,\gd_0)$ and $\tl \Gl_2=\tl \Gl_2(N,p,q,\Gw,\norm{u}_{L^1(D_{\frac{\gd_0}{2}})})$.
\es \medskip

\noindent{\bf Remark.} If $H$ satisfies \eqref{add} with $p \geq \frac{q}{2-q}$ then \eqref{est2-add} can be improved. Indeed, it follows from Keller-Osserman estimate that there exists a constant $C_{N,p}$ depending only on $N$ and $p$ such that
$$ u(x) \leq C_{N,p}\,\rho(x)^{-\gb_2} \forevery x \in \Gw. $$


Denote  by $G^\Gw$ (resp. $P^\Gw$)  the Green kernel (resp. the Poisson kernel) in $\Gw$, with corresponding operators $\BBG^\Gw$ (resp. $\BBP^\Gw$). We also denote by $\GTM_{\rho^\ga}(\Gw)$, $\ga \in [0,1]$, the space of Radon measures $\gm$ on $\Gw$ satisfying $\myint{\Gw}{}\rho^\ga(x) d|\gm|<\infty$, by $\GTM(\prt \Gw)$ the space of bounded Radon measures on $\prt \Gw$ and by $\GTM_+(\prt \Gw)$ the positive cone of $\GTM(\prt \Gw)$.

Denote  $L^p_w(\Gw;\tau)$, $1 \leq p < \infty$, $\tau \in \GTM_+(\Gw)$, the weak $L^p$ space defined as follows: a measureable function $f$ in $\Gw$ belongs to this space if there exists a constant $c$ such that
\bel{distri} \gl_f(a;\tau):=\tau(\{x \in \Gw: |f(x)|>a\}) \leq ca^{-p}, \forevery a>0. \ee
The function $\gl_f$ is called the distribution function of $f$ (relative to $\tau$). For $p \geq 1$, denote
$$ L^p_w(\Gw;\tau)=\{ f \text{ Borel measurable}: \sup_{a>0}a^p\gl_f(a;\tau)<\infty\} $$
and
\bel{semi} \norm{f}^*_{L^p_w(\Gw;\tau)}=(\sup_{a>0}a^p\gl_f(a;\tau))^{\frac{1}{p}}. \ee
The $\norm{.}_{L^p_w(\Gw;\tau)}$ is not a norm, but for $p>1$, it is equivalent to the norm
\bel{normLw} \norm{f}_{L^p_w(\Gw;\tau)}=\sup\left\{ \frac{\int_{\gw}|f|d\tau}{\tau(\gw)^{1/p'}}:\gw \sbs \Gw, \gw \text{ measurable }, 0<\tau(\gw)<\infty \right\}. \ee
More precisely,
\bel{equinorm} \norm{f}^*_{L^p_w(\Gw;\tau)} \leq \norm{f}_{L^p_w(\Gw;\tau)} \leq \myfrac{p}{p-1}\norm{f}^*_{L^p_w(\Gw;\tau)} \ee
The following usefull estimates involving Green and Poisson operators can be found in \cite{BVi} (see also \cite{MVbook}, \cite{V1} and \cite{V2}).
\bprop{P1} For any $\ga\in [0,1]$, there exist a positive constant $c_1$ depending on $\ga$,  $\Gw$ and $N$ such that
\bel{E1-1} \BA{lll}
\norm{\BBG^{\Gw}[\gn]}_{L^1(\Gw)}+\norm{\BBG^{\Gw}[\gn]}_{L_w^{\frac{N+\ga}{N+\ga-2}}(\Gw;\rho^\ga dx)}+\norm{\nabla\BBG^{\Gw}[\gn]}_{L_w^{\frac{N+\ga}{N+\ga-1}}(\Gw;\rho^\ga dx)} \\ [3mm]
\phantom{---------------------}
\leq c_1\norm \gn_{\mathfrak M_{\rho^\ga}(\Gw)},
\EA \ee
\bel{E1-3} \BA{lll}
\norm{\BBP^{\Gw}[\gm]}_{L^1(\Gw)}+\norm{\BBP^{\Gw}[\gm]}_{L_w^{\frac{N}{N-1}}(\Gw)}+\norm{\nabla\BBP^{\Gw}[\gm]}_{L_w^{\frac{N+1}{N}}(\Gw;\rho dx)}\leq c_1\norm \gm_{\GTM(\prt \Gw)},
\EA \ee
for any $\gn\in \GTM_{\rho^\ga}(\Gw)$ and any $\gm\in \GTM(\prt\Gw)$ where
$$ \norm \gn_{\mathfrak M_{\rho\ga}(\Gw)}:=\myint{\Gw}{}\rho^\ga(x) d|\gn| \q \text{and} \q \norm \gm_{\GTM(\prt \Gw)}=\myint{\prt \Gw}{}d|\gm|. $$
\es\medskip

\section{Boundary value problem with measures and boundary trace}
\subsection{The Dirichlet problem}
We first prove a regularity result in the subcritical case.
\blemma{reg-subcritical} Assume either $H$ satisfies \eqref{multi} with $0<N(p+q-1)<p+1$ or $H$ satisfies \eqref{add} with $m_{p,q}<p_c$. Let $u$ be a positive solution of \eqref{A0} (in the sense of distribution) satisfying $u \in L^s_{loc}(\Gw)$ for any $1<s<\frac{N+1}{N-1}$ and $\abs{\nabla u} \in L^r_{loc}(\Gw)$ for any $1<r<\frac{N+1}{N}$. Then $u \in C^2(\Gw)$. \es
\Proof We provide here the proof in the case $H$ satisfies \eqref{multi}. The case $H$ satisfies \eqref{add} follows by some modifications.The proof is based on bootstrap argument. We put $f=- u^p\abs{\nabla u}^q$. Let $s_1>1$ ($s_1$ will be determined later on) and $K \sbs \sbs \Gw$. By Holder inequality, for $r_1>1$ (will be made precise later),
\bel{hol} \myint{K}{}\abs{f}^{s_1}dx \leq \left( \myint{K}{}u^{ps_1r_1}dx \right)^{\frac{1}{r_1}}\left(\myint{K}{} \abs{\nabla u}^{qs_1r_1'}dx \right)^{\frac{1}{r_1'}}. \ee
We will choose $s_1$ and $r_1$ such that $ (N-1)ps_1r_1 \leq N+1$ and $Nqs_1r_1' \leq N+1$.
It is sufficient to choose $s_1$ and $r_1$ such that
$$ 1<s_1<\myfrac{N+1}{p(N-1)+qN} \quad \text{and} \quad \myfrac{N+1}{N+1-qNs_1}<r_1<\myfrac{N+1}{p(N-1)s_1}. $$
Since $u \in W^{1,s}_{loc}(\Gw)$ for any $1<s<\frac{N+1}{N}$, by Sobolev imbedding, $u \in L^{s_1}_{loc}(\Gw)$. It follows from interior regularity result for elliptic equations that $u \in W^{2,s_1}_{loc}(\Gw)$. Again, by Sobolev imbedding,
$ u \in L^{\frac{Ns_1}{N-2s_1}}_{loc}(\Gw)$ and $\abs{\nabla u} \in L^{\frac{Ns_1}{N-s_1}}_{loc}(\Gw)$ if $ s_1<\frac{N}{2}$.

Next, let $s_2>1$ (will be determined later on). By Holder inequality, for $r_2>1$ (will be made precise later), \eqref{hol} remains true with $s_1$ and $r_1$ replaced by $s_2$ and $r_2$ respectively. We will choose $s_2$ and $r_2$ such that
$ (N-2s_1)ps_2r \leq Ns_1$ and $(N-s_1)qs_2r' \leq Ns_1$. It is sufficient to choose $s_2$ and $r_2$ such that
$$ 1<\myfrac{(N+1)s_1}{p(N-1)+qN}<s_2<\myfrac{Ns_1}{p(N-2s_1)+q(N-s_1)} $$
$$ \myfrac{Ns_1}{Ns_1-q(N-s_1)s_2}<r_2<\myfrac{Ns_1}{p(N-2s_1)s_2}. $$
Then
$$ s_2 -s_2> \myfrac{N+1-p(N-1)-qN}{p(N-1)+qN}>0. $$
We can choose $s_2$ close  $\frac{(N+1)s_1}{p(N-1)+qN}$ enough that $u \in L^{s_2}_{loc}(\Gw)$; hence $u \in W^{2,s_2}_{loc}(\Gw)$. By Sobolev imbedding, $ u \in L^{\frac{Ns_1}{N-2s_2}}_{loc}(\Gw)$ and $\abs{\nabla u} \in L^{\frac{Ns_1}{N-s_2}}_{loc}(\Gw)$ if  $s_2<\frac{N}{2}$.
Next by iterating the process, we can define a sequence $\{s_k\}$ such that $u \in W^{2,s_k}_{loc}(\Gw)$ and
$$ s_k>1+k\myfrac{N+1-p(N-1)-qN}{p(N-1)+qN}. $$
Hence we can find $k_0$ large enough such that $s_{k_0}>N$ and $u \in W^{2,s_{k_0}}_{loc}(\Gw)$. By Sobolev imbedding, $u \in C^2(\Gw)$. \qed \medskip

We now turn to the \medskip

\noindent {\bf Proof of Theorem A.} We deal with the case when $H$ satisfies \eqref{multi}. The case $H$ satisfies \eqref{add} is simpler and can be treated in a similar way.  \smallskip

Let $\{\gm_n\}$ be a sequence of  positive functions in $C^1(\prt \Gw)$ such that $\{\gm_n\}$ converges to $\gm$ in the weak sense of measures and $\norm{\gm_n}_{L^1(\prt \Gw)} \leq c_2\norm{\gm}_{\GTM(\prt \Gw)}$ for all $n$, where $c_2$ is a positive constant independent of $n$. Consider the following problem
	\bel{v_n} \left\{ \BA{lll}
		- \Gd v +  (v+ \BBP^\Gw[\gm_n])^p|\nabla (v + \BBP^{\Gw}[\gm_n])|^q = 0 \qq &\text{ in } \Gw\\
		\phantom{- \Gd  +g(\abs{\nabla (v + \BBP^{\Gw}[\gm_n])}),,,----- }
		v = 0  &\text{ on } \prt \Gw.
	\EA \right. \ee
It is easy to see that $0$ and $-\BBP^{\Gw}[\gm_n]$ are respectively supersolution and subsolution of \eqref{v_n}. By \cite[Theorem 6.5]{KaKr} there exists a solution $v_n \in W^{2,p}(\Gw)$ with $1<p<\ity$ to problem \eqref{v_n} satisfying $-\BBP^{\Gw}[\gm_n] \leq v_n \leq 0$. Thus $u_n = v_n + \BBP^{\Gw}[\gm_n]$ is a solution of
	\bel{u_n} \left\{ \BA{lll}
		- \Gd u_n +   u_n^p\abs{\nabla u_n}^q = 0 \qq &\text{ in } \Gw\\
		\phantom{- \Gd  + g(\abs{\nabla u_n}),,- }
		u_n = \gm_n  &\text{ on } \prt \Gw.
	\EA \right. \ee
By the maximum principle, such solution is the unique solution of \eqref{u_n}. \medskip

\noindent{\bf Assertion 1:} $\{u_n\}$ and $\{\abs{\nabla u_n}\}$ remain uniformly bounded respectively in $ L_w^\frac{N}{N-1}(\Gw)$ and $L_w^\frac{N+1}{N}(\Gw;\rho dx)$.

Let $\gx$ be the solution to
	 \bel{eta} - \Gd \gx = 1 \text{ in } \Gw, \q \gx = 0  \text{ on } \prt \Gw, \ee
then there exists a constant $c_3>0$ such that $c_3^{-1}<-\frac{\prt \gx}{\prt \bf n} < c_3$ on $\prt \Gw$ and $c_3^{-1}\rho \leq \gx \leq c_3\rho$ in $\Gw$.
By multiplying the equation in \eqref{u_n} by $\gx$ and integrating on $\Gw$, we obtain
	\bel{EM2-1} \myint{\Gw}{}u_n dx +  \myint{\Gw}{}u_n^p\abs{\nabla u_n}^q\rho dx \leq c_4\norm{\gm}_{\GTM(\prt \Gw)} \ee
where $c_4$ is a positive constant independent of $n$. From \rprop{P1} and by noticing that $u_n \leq \BBP^{\Gw}[\gm_n]$, we get
	\bel{EM2-2} \norm{u_n}_{L_w^\frac{N}{N-1}(\Gw)} \leq \norm{\BBP^\Gw[\gm_n]}_{L_w^{\frac{N}{N-1}}(\Gw)} \leq c_1\norm{\gm_n}_{L^1(\prt \Gw)} \leq c_1c_2\norm{\gm}_{\GTM(\prt \Gw)}. \ee
Set $f_n= -   u_n^p\abs{\nabla u_n}^q$ then $f_n \in L_\rho^1(\Gw)$ and $u_n$ satisfies $-\Gd u_n = f_n$ in $\Gw$, $u_n=0$ on $\prt \Gw$. Again, from \rprop{P1} and  \eqref{EM2-1}, we derive that
	\bel{EM2-3} \norm{\nabla u_n}_{L_w^\frac{N+1}{N}(\Gw,\rho dx)} \leq c_1\left(\norm{f_n}_{L_\rho^1(\Gw)} + \norm{\gm_n}_{L^1(\prt \Gw)}\right)  \leq c_4'\norm{\gm}_{\GTM(\prt \Gw)} \ee
where $c_4'$ is a positive constant depending only on $\Gw$ and $N$.
Thus the assertion 1 follows from \eqref{EM2-2} and \eqref{EM2-3}. \medskip
	
By regularity results for elliptic equations \cite{Mi}, there exist a subsequence, still denoted by $\{u_n\}$, and a function $u$ such that $\{u_n\}$ and $\{|\nabla u_n|\}$ converges to $u$ and $|\nabla u|$  a.e. in $\Gw$.\medskip

\noindent {\bf Assertion 2:} $\{u_n\}$ converges to $u$ in $L^1(\Gw)$.

Indeed, since $u_n \in L_w^{\frac{N}{N-1}}(\Gw)$, it follow from \eqref{normLw} that if $G\subset\Gw$ be a Borel subset, then
\bel{EM2-5} \myint{G}{}u_ndx\leq |G|^{\frac{1}{N}}\norm{u_n}_{L_w^{\frac{N}{N-1}}(\Gw)}
\leq c_1c_2|G|^{\frac{1}{N}}\norm\gm_{\GTM(\prt\Gw)}.
\ee
Hence $\{u_n\}$ is uniformly integrable. Therefore the assertion 2 follows from Vitali's convergence theorem.\medskip

\noindent {\bf  Assertion 3:} $\{u_n^p|\nabla u_n|^q\}$ converges to $u^p|\nabla u|^q$ in $L_\rho^1(\Gw)$.

Indeed, let $G$ be a Borel in $\Gw$, $\ell>0$, $\gl>0$ and write
\bel{S1} \myint{G}{}u_n^p|\nabla u_n|^q\rho dx= I_1+I_2+I_3+I_4 \ee
where
\bel{S2} \BA{lll}
I_1:=\myint{G\cap\{x:u_n\leq\ell, |\nabla u_n(x)|\leq\gl\}}{}u_n^p|\nabla u_n|^q \rho \,dx,  \\[3mm] I_2:=\myint{G\cap\{x:u_n>\ell, |\nabla u_n(x)|\leq\gl\}}{}u_n^p|\nabla u_n|^q \rho \, dx  \\[3mm]
I_3:=\myint{G\cap\{x:u_n\leq\ell, |\nabla u_n(x)|>\gl\}}{}u_n^p|\nabla u_n|^q \rho \, dx,  \\[3mm]
I_4:=\myint{G\cap\{x:u_n>\ell, |\nabla u_n(x)|>\gl\}}{}u_n^p|\nabla u_n|^q \rho \, dx.
\EA \ee
We first notice that
\bel{EQ0} I_1 \leq \ell^p \gl^q \myint{G}{}\rho dx. \ee
Next put $A_n(t)=\{x \in \Gw: u_n>t\}$, $t>0$ and $a_n(t)=\myint{A_n(t)}{}\rho dx$. Since $0 \leq u_n \leq \BBP^\Gw[\gm_n]$ and $\BBP^\Gw[\gm_n] \in L_w^{\frac{N+1}{N-1}}(\Gw;\rho dx)$, it follows that $u_n \in L_w^{\frac{N+1}{N-1}}(\Gw;\rho dx)$ . By \eqref{semi} and \eqref{equinorm}, we get
$$a_n(t) \leq t^{-\frac{N+1}{N-1}}(\norm{u_n}^*_{L_w^{\frac{N+1}{N-1}}(\Gw;\rho dx)})^{\frac{N+1}{N-1}} \leq  t^{-\frac{N+1}{N-1}}(\norm{u_n}_{L_w^{\frac{N+1}{N-1}}(\Gw;\rho dx)})^{\frac{N+1}{N-1}}.  $$
Combining the above inequality with \eqref{E1-3} yields
\bel{an} a_n(t)  \leq t^{-\frac{N+1}{N-1}}(c_1c_2\norm{\gm}_{\GTM(\prt \Gw)})^{\frac{N+1}{N-1}}:= c_5t^{-\frac{N+1}{N-1}}. \ee
By integration by part and \eqref{an}, we obtain
$$ \BA{lll} \myint{\{x \in \Gw: u_n>\ell\}}{}u_n^p\,\rho \,dx &= -\myint{\ell}{\infty}t^pda_n(t)=\ell^p a_n(\ell) + p\myint{\ell}{\infty}a_n(t)t^{p-1}dt \\
&\leq \ell^pa_n(\ell) + c_5p\myint{\ell}{\infty}t^{p-1-\frac{N+1}{N-1}}dt \\
&\leq \myfrac{c_5(N+1)}{N+1-p(N-1)}\ell^{p-\frac{N+1}{N-1}}
\EA $$
Hence
$$I_2 \leq \myfrac{c_5(N+1)}{N+1-p(N-1)}\gl^q\ell^{p-\frac{N+1}{N-1}}.$$
If $\gl=\ell^{\frac{N}{N-1}}$ then
\bel{EQ1} I_2 \leq \myfrac{c_5(N+1)}{N+1-p(N-1)}\ell^{\frac{p(N-1)+qN-(N+1)}{N-1}}. \ee
Similarly, we get
\bel{EQ2} I_3 \leq \myfrac{c_6(N+1)}{N+1-qN}\ell^{\frac{p(N-1)+qN-(N+1)}{N-1}}. \ee
Fix $r$ such that
\bel{EQ3} \myfrac{N+1}{N+1-qN}<r<\myfrac{N+1}{p(N-1)} \ee
then by Holder inequality,
$$  I_4 \leq \left(\myint{\{x: u_n>\ell\}}{}u^{pr}dx\right)^{\frac{1}{r}}\left(\myint{\{x:|\nabla u_n|>\gl\}}{}\abs{\nabla u}^{qr'}dx\right)^{\frac{1}{r'}}. $$
Due to the choice of $r$ in \eqref{EQ3}, $r'<\frac{N+1}{Nq}$ where $\frac{1}{r}+\frac{1}{r'}=1$. Therefore
\bel{EQ4} I_4 \leq c_7\ell^{p-\frac{N+1}{(N-1)r}}\gl^{q-\frac{N+1}{Nr'}}=c_7\ell^{\frac{p(N-1)+qN-(N+1)}{N-1}} \ee
where $ c_7=c_5(N+1)(N+1-pr(N-1))^{-r}(N+1-qr'N)^{-r'}$.
Combining \eqref{S1}, \eqref{EQ0}-\eqref{EQ4} yields
\bel{EQ5}  \myint{G}{}u_n^p|\nabla u_n|^q\rho dx \leq \ell^{\frac{p(N-1)+qN}{N-1}}\myint{G}{}\rho dx + c_8 \ell^{\frac{p(N-1)+qN-(N+1)}{N-1}}. \ee
For any $\vge>0$, since $p(N-1)+qN<N+1$, we fix $\ell$ large enough that the second term on the right-hand side of \eqref{EQ5} is smaller than $\frac{\vge}{2 }$. Therefore,
	$$ \myint{G}{}\rho \, dx <\frac{\vge}{2 }\ell^{-\frac{p(N-1)+qN}{N-1}} \Lra   \myint{G}{}u_n^p|\nabla u_n|^q\rho dx < \vge. $$
Thus the assertion 3 is a consequence of Vitali's convergence theorem. \medskip

For every $\zeta \in C_0^2(\ovl \Gw)$, we have
	\bel{EM2-8} \myint{\Gw}{}(-u_n\Gd \zeta +  u_n^p \abs{\nabla u_n}^q\zeta)dx=-\myint{\prt \Gw}{}\gm_n\myfrac{\prt \zeta}{\prt \bf n}dS.
	\ee
From assertion 2 and assertion 3, by taking into account that $|\zeta|\leq c\rho$ in $\Gw$, we can pass to the limit in each term in \eqref{EM2-8} and obtain \eqref{A2}; so $u$ is a solution of \eqref{A1}.By \rprop{P1} $u \in L_w^\frac{N}{N-1}(\Gw)$ and  $\abs{\nabla u} \in L_w^\frac{N+1}{N}(\Gw;\rho dx)$. \medskip

Next, let $\{\gm_n\}$ be a sequence of positive bounded measures on $\prt\Gw$ which converges to a positive bounded $\gm$ in the weak sense of measures and $\{u_{\gm_n}\}$ is be a sequence of corresponding solutions of \eqref {u_n}. Then by using the same argument as in assertion 2 and assertion 3, we deduce that there exists a subsequence such that $\{ u_{\gm_{n_k}}\}$ converges to a solution $ u_\gm$ of \eqref {A1} in $L^1(\Gw)$ and $\{H(x,u_{\gm_{n_k}},\nabla u_{\gm_{n_k}})\}$ converges to $H(x,u,\nabla u)$ in $L^1_\rho(\Gw)$. \qed \medskip

A variant of the stability result in theorem A is the following
\bcor{var-stab}
Assume either $H$ satisfies \eqref{multi} with $0<N(p+q-1)<p+1$ or $H$ satisfies \eqref{add} with $m_{p,q}<p_c$. Let $\{\gd_n\}$ be a decreasing sequence converging to $0$, $\gm$ is a bounded positive measure on $\prt \Gw$ and $\{\gm_n\}$ is a sequence of bounded positive measure on $\Gs_{\gd_n}$ converging to $\gm$ in the weak sense of measures and $\{u_{\gm_n}\}$ be a sequence of corresponding solutions of \eqref {u_n} in $D_{\gd_n}$. Then there exists a subsequence such that $\{ u_{\gm_{n_k}}\}$ converges in $L^1(\Gw)$ to a solution $ u_\gm$ of \eqref {A1} and $\{H(x,u_{\gm_{n_k}},\nabla u_{\gm_{n_k}})\}$ converges to $H(x,u,\nabla u)$ in $L^1_\rho(\Gw)$.
\es
\Proof As above, we consider the case $H$ satisfies \eqref{multi} because the case $H$ satisfies \eqref{add} follows by similar argument. We extend $ u_{\gm_n}$ and $|\nabla u_{\gm_n}|$ by zero outside $\ovl D_{\gd_n}$ and still denote them by the same expressions. By regularity results for elliptic equations \cite{Mi}, there exist a subsequence, still denoted by $\{u_{\gn_n}\}$, and a function $u$ such that $\{u_{\gm_n}\}$ and $\{|\nabla u_{\gm_n}|\}$ converges to $u$ and $|\nabla u|$  a.e. in $\Gw$. Let $G\subset\Gw$ be a Borel set and put $G_n=G\cap D_{\gd_n}$. By using similar argument as in assertion 2 in the proof of theorem A, thanks to the estimate $||\BBP^{\Gw}[\gm]|_{_{\Gs_{\gd_n}}}||_{L^1(\Gs_{\gd_n})}\leq c_7\norm{\gm}_{\mathfrak M(\Gs)}$, we derive
\bel{EM2-5+1}\BA {l} \myint{G_n}{}u_{\gm_n} dx\leq |G_n|^{\frac{1}{N}}\norm{u_{\gm_n}}_{L_w^{\frac{N}{N-1}}(D_{\gd_n})}
\leq c_1c_2|G_n|^{\frac{1}{N}}\norm{\BBP^\Gw[\gm]|_{_{\Gs_{\gd_n}}}}_{L^1(\Gs_{\gd_n})}\\[4mm]\phantom{\myint{G_n}{}u_\gd dx}
\leq c_1c_2c_7|G|^{\frac{1}{N}}\norm{\gm}_{\GTM(\Gs)}.
\EA\ee
Hence, $\{u_{\gm_n}\}$ is uniformly integrable. Therefore, due to Vitali's convergence theorem, up to a subsequence,$\{u_{\gm_n}\}$ converges to $u$ in $L^1(\Gw)$.

Set $\rho_n(x):=(\rho(x)-\gd_n)_+$. By using a similar argument as in Assertion 3 of the proof of Theorem A and taking into account that $\int_{G_n}^{}\rho_n dx \leq \int_{G}^{}\rho dx$, we obtain that for any $\vge>0$ there exists $\ell>0$ large enough, independent of $n$, such that
\bel{EM2-6} \myint{G_n}{}\rho_n dx <\frac{\vge}{2 }\ell^{-\frac{p(N-1)+qN}{N-1}} \Lra  \myint{G_n}{}u_{\gm_n}^p|\nabla u_{\gm_n}|^q\rho_n dx < \vge. \ee
Therefore, by Vitali's convergence, up to a subsequence, $\{ u_{\gm_n}^p |\nabla u_{\gm_n}|^q\}$ converges to $u^p|\nabla u|^q$ in $L^1_\rho(\Gw)$.

Finally, if $\gz\in C_0^2(\Gw)$ we denote by $\gz_n$ the solution of
\bel{EM2-9} -\Gd \gz_n=-\Gd\gz \text{ in } D_{\gd_n}, \q \gz_n=0 \text{ on } \prt D_{\gd_n}. \ee
Then $\gz_n \in C_0^2(\ovl \Gw_{\gd_n})$, $\gz_n \to \gz$ in $C^2(\Gw)$ and $\sup_n\norm{\gz_n}_{C^2(\ovl \Gw_{\gd_n})}<\infty$. Furthermore,
	\bel{EM2-10} \myint{D_{\gd_n}}{}(- u_{\gm_n}\Gd \zeta_n +   u_{\gm_n}^p\abs{\nabla u_{\gm_n}}^q\zeta_n)dx=-\myint{\Gs_{\gd_n}}{}\myfrac{\prt \zeta_n}{\prt \bf n}d \gm_n
	\ee
By letting $n \to \infty$, we deduce that $u$ is a solution of \eqref{A1}. \qed

\noindent {\bf Remark.}  Let $\gm \in \GTM_+(\prt \Gw)$ and $u$ is a positive solution of \eqref{A1}. It follows from \rprop{est-multi} that there exists a constant $c$ depending on $N$, $p$, $q$, $\Gw$ and $\norm{\gm}_{\GTM(\prt \Gw)}$ such that
	\bel{est1*} u(x) \leq c\rho(x)^{-\gb_1} \forevery x \in \Gw, \ee
	\bel{est2*}|\nabla u(x)| \leq c\rho(x)^{-\gb_1-1} \forevery x \in \Gw. \ee
\bdef{potential} A nonnegative superharmonic function is called a $\Gd$-potential if its largest harmonic minorant is zero.
\es
\bprop{nontan} Assume either $H$ satisfies \eqref{multi} with $0<N(p+q-1)<p+1$ or $H$ satisfies \eqref{add} with $m_{p,q}<p_c$. Let $\gm \in \GTM_+(\prt \Gw)$. If $u$ is a positive solution of \eqref{A1} then
	\bel{nontang} \lim_{x \to y}\myfrac{u(x)}{\BBP^\Gw[\gm](x)}=1 \q non-tangentially, \, \gm-a.e. \ee
\es
\Proof Put $v_\gm=\BBP^\Gw[\gm]-u$ then $v_\gm>0$ and $-\Gd v_\gm=u_\gm^p|\nabla u_\gm|^q \geq 0$ in $\Gw$. It means $v_\gm$ is a positive superharmonic function in $\Gw$. By Riesz Representation Theorem (see \cite{Ma}), $v_\gm$ can be written as follows: $ v_\gm= v_h + v_p $ where $v_h$ is a nonnegative harmonic function and $v_p$ is a $\Gd$-potential (see \cite{Ma} for more details). Since the boundary trace of $v_\gm$  is a zero measure, it follows the boundary trace of $v_h$ and $v_p$ is zero measure. Hence $v_h = 0$ in $\Gw$, therefore $v_\gm=v_p$. By \cite[Theorem 2.11 and Lemma 2.13]{Ma}, we derive \eqref{nontang}. \qed

\subsection{Boundary trace}
\bdef{moderate} A positive solution $u$ of \eqref{A0} is {\bf moderate} is $u \in L^1(\Gw)$ and $H(.,u,\nabla u) \in L^1_\rho(\Gw)$. \es

It's clear to see that
\bprop{moderate} The following statements are equivalent \medskip

\noindent i) $u$ is a moderate solution of \eqref{A0}. \medskip

\noindent ii) There exists $\gm \in \GTM_+(\prt \Gw)$ such that $u$ is a solution of \eqref{A1}. \medskip

\noindent iii) $u$ is bounded  from above by an harmonic function in $\Gw$.
\es

\bdef{tra} Assume $\gm \in \GTM(\Gs)$ and $\gm_\gd \in \GTM(\Gs_\gd)$ for each $\gd \in (0,\gd_0)$. We say that $\gm_\gd\to \gm$ as $\gd \to 0$ in the sense of weak convergence of measures if
\bel{E2}
\lim_{\gd\to 0}\myint{\Gs_\gd}{}\gf(\gs(x))d\gm_\gd=\myint{\Gs}{}\gf \,d\gm\qq\forall\gf\in C_c(\Gs).
\ee
A function $u\in C(\Gw)$ of \eqref{A0} possesses a measure boundary trace $\gm\in \mathfrak M(\Gs)$ if
\bel{E3}
\lim_{\gd\to 0}\myint{\Gs_\gd}{}\gf(\gs(x))u(x)dS=\myint{\Gs}{}\gf \,d\gm\qq\forall\gf\in C_c(\Gs).
\ee
Similarly, if $A$ is a relatively open subset of $\Gs$, we say that $u$ possesses a trace $\gm$ on $A$ in the sense of weak convergence of measures if $\gm \in \GTM(A)$ and \eqref{E3} holds for every $\gf \in C_c(A)$.
\es

By adapting the proof of \cite[Cor 2.3]{MV4} to \eqref{A0}, we obtain
\bprop{tra-reg} Let $u\in C^2(\Gw)$ be a positive solution of \eqref{A0}. Suppose that for some $z\in\prt\Gw$ there exists an open neighborhood $U$ such that
\bel{E4}
\myint{U\cap \Gw}{}H(x,u,\nabla u) \rho\,  dx<\infty.
\ee
Then $u\in L^1(K \cap \Gw)$ for every compact set $K\sbs U$ and there exists a positive Radon measure $\gn$ on $\Gs\cap U$ such that
\bel{E5}
\lim_{\gd\to 0}\myint{\Gs_\gd\cap U}{}\gf(\gs(x))u(x)dS=\myint{\Gs \cap U}{}\gf \,d\gn\qq\forall\gf\in C_c(\Gs\cap U).
\ee
\es

\bdef{reg-sing} Let $u\in C^2(\Gw)$ be a positive solution of \eqref{A0}. A point $z\in\prt\Gw$ is a regular boundary point of $u$ if there exists an open neighborhood $U$ of $z$ such that \eqref{E4} holds. The set of regular points is denoted by $\CR(u)$. Its complement $\CS(u)=\prt\Gw\setminus\CR(u)$ is called the singular boundary set of $u$.
\es

Clearly $\CR(u)$ is relatively open and there exists a positive Radon measure $\gm$ on $\CR(u)$ such that $u$ admits $\gm:=\gm(u)$ as a measure boundary trace on $\CR(u)$ and $\gm(u)$ is uniquely determined. The couple $(\CS(u),\gm)$ is called the {\it boundary trace} of $u$ and denoted by $tr_{\prt\Gw}(u)$. \medskip

Concerning $\CS(u)$, we get the following result by employing \cite[Lemma 2.8]{MV4}.
\bprop{L1}Let $u\in C^2(\Gw)$ be a positive solution of \eqref{A0} with the singular boundary set $\CS(u)$. If  $z\in
\CS(u)$ is such that there exists an open neighborhood $U$ of $z$ such that $u\in L^1(U\cap\Gw)$, then for every neighborhood $V$ of $z$ there holds
\bel{E6}
\lim_{\gd\to 0}\myint{\Gs_\gd\cap V}{}u \, dS=\infty.
\ee
\es
\bth{gen+M} Assume either $H$ satisfies \eqref{multi} with $0<N(p+q-1)<p+1$ or $H$ satisfies \eqref{add} with $m_{p,q}<p_c$. If $u\in C^2(\Gw)$ is a positive solution of \eqref{A0}, then \eqref{E6} holds for every $z\in \CS(u)$.\es
 \Proof By translation we assume $z=0\in\CS(u)$ and \eqref{E6} does not hold. We proceed by contradiction, assuming that there exists an open neighborhood $G$ of $0$ such that
  \bel{E17}
  \liminf_{\gd\to 0}\myint{\Gs_\gd\cap G}{}u\,dS<\infty.
\ee
It follows from \rprop{L1} that if $U$ is a neighborhood  of $z$ then $\int_{\Gw\cap U}^{}u\,dx=\infty$, which leads to $\limsup_{\gd\to 0}\int_{\Gs_\gd\cap U}^{}u\,dS=\infty$. For each $n\in\BBN_*$, we take $U=B_{\frac{1}{n}}(0)$. Then there exists a sequence $\{\gd_{n,m}\}_{m\in\BBN}$ tending to $0$ as $m \to \infty$ such that $  \lim_{m\to \infty}\int_{\Gs_{\gd_{n,m}}\cap B_{\frac{1}{n}}(0)}^{}u\,dS=\infty$.
Then, for any $k>0$, there exists $m_k:=m_{n,k}\in\BBN$ such that
 \bel{E20}
 m\geq m_k\Longrightarrow \myint{\Gs_{\gd_{n,m}}\cap B_{\frac{1}{n}}(0)}{}udS\geq k
 \ee
and $m_{n,k}\to\infty$ when $n\to \infty$.  In particular there exists $t:=t(n,k)>0$ such that
 \bel{E21}
 \myint{\Gs_{\gd_{n,m_k}}\cap B_{\frac{1}{n}}(0)}{}\inf\{u,t\}dS=k.
 \ee
By the comparison principle $u$ is bounded from below in $D_{\gd_{n,m_k}}$ by the solution $v:=v_{\gd_{n,m_k}}$ of
 \bel{E22}\left\{\BA {ll}
 -\Gd v+ H(x,v,\nabla v)=0\qq&\text{in }D_{\gd_{n,m_k}}\\[2mm]
 \phantom{ -\Gd +g(|\nabla v|),,,,}
v=\inf\{u,t\}\qq&\text{on }\Gs_{\gd_{n,m_k}}.
\EA\right. \ee
When $n\to\infty$, $\inf\{u,t(n,k)\}dS$ converges in the weak sense of measures to $k\gd_0$. By \rcor{var-stab} there exists a subsequence, still denoted by $\{v_{\gd_{n,m_k}}\}_n$,  such that $v_{\gd_{n,m_k}}\to u_{k,0}^\Gw$ when $n\to\infty$ where $u_{k,0}^\Gw$ is the unique solution of \eqref{APE0} and consequently
$u\geq u_{k,0}^\Gw$ in $\Gw$. Therefore,
 \bel{E23}
\liminf_{\gd\to 0}\myint{\Gs_{\gd}}{}u\gz \,dS\geq \lim_{\gd\to 0}\myint{\Gs_{\gd}}{}u_{k,0}^\Gw\gz \, dS=k \ee
for any nonnegative $\gz\in C^{\infty}(\BBR^N)$ such that $\gz=1$ in a neighborhood of $0$. Since $k$ is arbitrary we obtain
 \bel{E24}
\liminf_{\gd\to 0}\myint{\Gs_{\gd}}{}u\,\gz \, dS=\infty
\ee
which contradicts \eqref{E17}.\qed

\section{Isolated boundary singularities}
\subsection{Weak singularities}
Let us give useful estimates involving solutions with isolated boundary singlarities which play a key role in the proof of uniquess and removability results in the sequel.
\blemma{estfunct2} Assume $u\in C(\overline \Gw\setminus \{0\})\cap C^2(\Gw)$ is a nonnegative solution of \eqref{A0} in $\Gw$ which vanishes on $\prt\Gw\setminus \{0\}$. \medskip

\noindent i) Assume $H$ satisfies \eqref{multi}. Then
\bel{C8} u(x)\leq \Gl_1|x|^{-\gb_1}\qq\forall x\in\Gw, \ee
\bel{C8a} |\nabla u(x)| \leq \Gl_3\abs{x}^{-\gb_1-1} \forevery x \in \Gw, \ee
\bel{C8b}  u(x)| \leq \tl \Gl_3\rho(x)\abs{x}^{-\gb_1-1} \forevery x \in \Gw \ee
where $\Gl_1$ is defined in \eqref{Lambda1}, $\Gl_3=\Gl_3(N,p,q,\Gw)$ and $\tl \Gl_3=\tl \Gl_3(N,p,q,\Gw)$ . \medskip

\noindent ii) Assume $H$ satisfies \eqref{add}. Then
\bel{C8-add} u(x)\leq \Gl_2|x|^{-\gb_2}\qq\forall x\in\Gw, \ee
\bel{C8a-add} |\nabla u(x)| \leq \Gl_4\abs{x}^{-\gb_2-1} \forevery x \in \Gw, \ee
\bel{C8b-add}  u(x)| \leq \tl \Gl_4\rho(x)\abs{x}^{-\gb_2-1} \forevery x \in \Gw \ee
where $\Gl_2=\Gl_2(p,q)$, $\Gl_4=\Gl_4(N,p,q,\Gw)$ and $\tl \Gl_4=\tl \Gl_4(N,p,q,\Gw)$.
\es
\Proof The proof is an adaptation of \cite[Lemma 3.3 and Lemma 3.4]{NV}. We deal only with the case where $H$ satisfies \eqref{multi} since the case $H$ satisfies \eqref{add} can be treated in a similar way. For $\ge>0$, we set
$$P_\ge(r)=\left\{\BA{ll}
0&\text{if }r\leq\ge\\
\frac{-r^4}{2\ge^3}+\frac{3r^3}{\ge^2}-\frac{6r^2}{\ge}+5r-\frac{3\ge}{2}\q&\text{if }\ge<r<2\ge\\
r-\frac{3\ge}{2}&\text{if } r \geq 2\ge
\EA\right.
$$
and let $u_\ge$ be the extension of $P_\ge(u)$ by zero outside $\Gw$. There exists $R_0$ such that $\Gw\subset B_{R_0}$. Since $0\leq P'_\ge(r)\leq 1$ and $P_\ge$ is convex, $u_\ge \in C^2(\BBR^N)$ and it satisfies $-\Gd u_\ge+u_\ge^p|\nabla u_\ge|^q\leq 0$. Furthermore $u_\ge$ vanishes in $ B^c_{R_0}$. For $R\geq R_0$ we set
$$U_{\ge,R}(x)=\Gl_1\left((|x|-\ge)^{-\gb_1}-(R-\ge)^{-\gb_1}\right)\qq\forall x\in B_R\setminus B_\ge,
$$
 then $-\Gd U_{\ge,R}+U_{\ge,R}^p|\nabla U_{\ge,R}|^q\geq 0$. Since $u_\ge$ vanishes on $\prt B_R$ and is finite on $\prt B_\ge$ it follows $u_\ge\leq U_{\ge,R}$ in $B_R \sms {\ovl B_\ge}$. Letting successively $\ge\to 0$ and $R\to\infty$ yields to \eqref{C8}.

For $\ell>0$, define $T^1_\ell[u](x)=\ell^{\gb_1}u(\ell x)$, $ x \in \Gw^\ell:=\ell^{-1}\Gw$. If $x_0\in\Gw$, we set $R=|x_0|$ and $u_R(x)=T^1_R[u](x)$. Then $u_R$ satisfies \eqref{A0} in $\Gw^{R}$. By \eqref{C8}, $\max\{|u_R(x)|:\frac{1}{2}\leq |x|\leq \frac{3}{2} \} \leq 2^{\gb_1}\Gl_1$. By \cite[Theorem 1]{Lib}, there exists $\Gl_3=\Gl_3(N,\Gw,p,q)$ such that $\max\{|\nabla u_R(x)|:\frac{3}{4}\leq |x|\leq \frac{5}{4} \} \leq \Gl_3$. In particular, $|\nabla u_R(x)| \leq \Gl_3$ with $|x|=1$. Hence $|\nabla u(x_0)| \leq \Gl_3|x_0|^{-\gb_1-1}$.

Finally, \eqref{C8b} follows from \eqref{C8} and \eqref{C8a}. \qed \medskip

Uniqueness can be obtained if $\gm$ is a bounded measure concentrated at a point on $\prt \Gw$.
\bth{unique} Assume either $H$ satisfies \eqref{multi} with $0<N(p+q-1)<p+1$ or $H$ satisfies \eqref{add} with $m_{p,q}<p_c$. Then for every $k >0$, there exists a unique solution,  denoted by $u^\Gw_{k,0}$, of the problem
           \bel{APE0} \left\{ \BA{lll}
	- \Gd u + H(x,u,\nabla u) &= 0 \qq &\text{in } \Gw\\
	\phantom{------,,,,,}
	u &= k\gd_0 &\text{on } \prt \Gw.
	\EA \right. \ee
Moreover, $u^\Gw_{k,0}(x)=k(1+o(1))P^\Gw(x,0)$ as $x \to 0$. Consequently the mapping $k \mapsto u^\Gw_{k,0}$ is increasing.
\es

The existence of a solution to \eqref{APE0} is guaranteed by Theorem A. The uniqueness is obtained due to the following lemma.

\blemma{APE-grad} Under the assumption of \rth{unique}, let $u$ be a solution to \eqref{APE0}. Then
\bel{AB4} \BBG^{\Gw}[H(.,u,\nabla u)](x)= o(P^\Gw(x,0)) \quad \text{as } x \to 0. \ee
\es
\Proof We prove \eqref{AB4} in the case $H$ satisfies \eqref{multi}. The case $H$ satisfies \eqref{add} can be treated in a similar way. \smallskip

\noindent {\bf Assertion:} For every $x \in \Gw$, there hold
\bel{APE-grad1} \abs{\nabla u(x)} \leq  \Gl_5 k\abs{x}^{-N} \ee
where $\Gl_5$ is a positive constant depending on $N,p,q,\Gw$.

Indeed, since $u$ is a solution of \eqref{APE0}, it follows from the maximum principle that $u \leq kP^\Gw(.,0) \leq kc_N|x|^{1-N}$ in $\Gw$ where $c_N$ is a positive constant depending on $N$ and $\Gw$. By adapting argument in \rlemma{estfunct2}, we obtain \eqref{APE-grad1}. \medskip

Next, it follows from \rlemma{APE-grad} that
\bel{G1} \BBG^\Gw[H(.,u,\nabla u)](x) \leq  c_8\myint{\Gw}{}G^\Gw(x,y)\abs{y}^{-(N-1)p-Nq}dy \forevery x \in \Gw. \ee
Since  $G^\Gw(x,y) \leq c_9\rho(x)\abs{x-y}^{-N}\min\{\abs{x-y},\rho(y)\}$ for every $x,y \in \Gw, x \neq y$  where $c_9=c_9(N,\Gw)$ (see \cite{MVbook}), we deduce, for $\ge_0 \in (0,1)$, that
	\bel{AB5} \BA{lll} G^\Gw(x,y) \leq c_9\rho(x)\rho(y)^{1-\ge_0}\abs{x-y}^{\ge_0-N} \forevery x,y\in \Gw, x \neq y,
	\EA \ee
The above estimate, joint with \eqref{G1}, implies that
	\bel{AB6} \BA{lll}\BBG^{\Gw}[H(.,u,\nabla u)](x) \\ \phantom{----}
\leq c_{10}\abs{x}^NP^\Gw(x,0)\myint{\BBR^N}{}\abs{x-y}^{\ge_0-N}\abs{y}^{1-(N-1)p-Nq-\ge_0}dy. \EA \ee
We fix $\vge_0$ such that $0<\vge_0<\min\left\{1,N+1-(N-1)p-Nq\right\}$. By the following identity (see \cite{LiLo}),
	\bel{AB7} \myint{\BBR^N}{}\abs{x-y}^{\vge_0-N}\abs{y}^{1-(N-1)p-Nq-\vge_0}dy = c_{11}\abs{x}^{N+1-(N-1)p-Nq} \ee
where $c_{11}=c_{11}(N,\vge_0)$, we obtain
$$\BBG^{\Gw}[H(.,u,\nabla u)](x) \leq \ga_1c_{10} c_{11}\abs{x}^{N+1-(N-1)p-Nq}P^\Gw(x,0).$$
Since $N+1-(N-1)p-Nq>0$, by letting $x \to 0$, we obtain \eqref{AB4}. \qed \medskip

\noindent{\bf Proof of \rth{unique}.} Let $u$ be a solution of \eqref{APE0} then $ u(x)=k\,P^\Gw(x,0)-\BBG^{\Gw}[H(.,u,\nabla u)](x)$.  From \eqref{AB4}, we obtain
	\bel{AB9} u(x)=k(1+o(1))P^\Gw(x,0) \quad \text{as } x \to 0, \ee
which, along with the comparison principle, implies the uniqueness of $u^\Gw_{k,0}$ and the monotonicity of $k\mapsto u^\Gw_{k,0}$. \qed \medskip

\bprop{compare} Assume either $H$ satisfies \eqref{multi} with $0<N(p+q-1)<p+1$ or $H$ satisfies \eqref{add} with $m_{p,q}<p_c$. Then for every $k>0$, there exists a positive constant $d_k$  depending on $N$, $p$, $q$, $k$ and $\Gw$ such that
	\bel{compare} d_k P^\Gw(x,0)<u^\Gw_{k,0}(x) < k P^\Gw(x,0) \forevery x \in \Gw. \ee
\es
\Proof The second inequality follows straightforward from comparison principle. In order to prove the first inequality, put $\CA=\{d>0: d\, P^\Gw(.,0)<u^\Gw_{k,0}\}$. Suppose by contradiction that $\CA=\emptyset$. Then for each $n \in \BBN$, there exists a point $x_n \in \Gw$ such that
\bel{contra}  n\, u^\Gw_{k,0}(x_n) <  P^\Gw(x_n,0). \ee
We may assume that $\{x_n\}$ converges to a point $x^* \in \ovl \Gw$. We deduce from \eqref{contra} that $x^* \notin \Gw$. Thus $x^* \in \prt \Gw$. By \rth{unique}, $x^* \in \prt \Gw \sms B_{\ge}(0)$ for some $\ge>0$. Following the notations in Section 2, denote by $\gs(x_n)$ the projection of $x_n$ on $\prt \Gw$. It follows from \eqref{contra} that
	$$  \myfrac{u^\Gw_{k,0}(\gs(x_n))-u^\Gw_{k,0}(x_n)}{\rho(x_n)} > \myfrac{1}{n}\myfrac{P^\Gw(\gs(x_n),0)- P^\Gw(x_n,0)}{\rho(x_n)}. $$
By letting $n \to \infty$, we obtain $ \frac{\prt u^\Gw_{k,0}}{\prt \bf n}(x^*) \geq 0$
which contradicts Hopf lemma. Thus $\CA \neq \emptyset$. Put $d_k=\max\CA$. By combining \eqref{C8} and boundary Harnack inequality, we deduce that $d_k$ depends on $N$, $p$, $q$, $k$ and $\Gw$.  \qed \medskip

\noindent {\bf Proof of Theorem B.} The proof follows from \rth{unique} and \rprop{compare}. \qed \medskip

The next result give us existence and uniqueness of weakly singular solution in the case that {\it $\Gw$ is unbounded domain}.
\bth{unbd} Assume either $H$ satisfies \eqref{multi} with $0<N(p+q-1)<p+1$ or $H$ satisfies \eqref{add} with $m_{p,q}<p_c$. Let  either $\Gw=\BBR_+^N:=\{x=(x',x_N):x_N>0\}$ or $\prt\Gw$ be compact with $0\in\prt\Gw$ ($\Gw$ is possibly unbounded). Then there exists a unique solution $u^\Gw_{k,0}$ to problem \eqref{APE0}.
\es
\Proof If $\prt\Gw$ is compact, for each $n \in \BBN$ large enough, $\prt \Gw \sbs B_n(0)$. We set $\Gw_n=\Gw\cap B_n(0)$ and denote by $u^{\Gw_n}_{k,0}$ the unique solution of
 	\bel {Pc+1}\left\{ \BA{lll} -\Gd u + H(x,u,\nabla u)  = 0 \qq &\text{in } \Gw_n \\
 	\phantom{ -\Gd  + u^p\abs{\nabla u}^q,,,}
 	u = k\gd_0 &\text{on } \prt \Gw_n.
 	\EA \right. \ee
Then by the maximum principle,
	 \bel {Pc+2}u^{\Gw_n}_{k,0}(x)\leq kP^{\Gw_n}(x,0)\qq\forall x\in \Gw_n.
	 \ee	
Thus $\{u^{\Gw_n}_{k,0}\}$ increase to a function $u^*$ which satisfies
	 	\bel {Pc+3}
	u^*(x)\leq kP^{\Gw}(x,0)\qq\forall x\in \Gw.
	\ee
By regularity theory, $\{\nabla u^{\Gw_n}_{k,0}\}_n$ converges locally uniformly in $\overline \Gw\setminus B_\ge(0)$ for any $\ge>0$ when $n\to\infty$, and  thus $u^*\in C(\overline\Gw\setminus\{0\})$ is a positive solution of  \eqref{APE0} in $\Gw$ vanishing on $\prt\Gw\setminus\{0\}$. The estimate \eqref{Pc+3} implies that the boundary trace of $u^*$ is a Dirac measure at $0$, which is in fact $k\gd_0$ due to  \eqref{AB9} for $\Gw_n$, \eqref{Pc+2} and \eqref{Pc+3}. Uniqueness also follows from these estimates.\qed

\subsection{Strong singularities} 
\noindent {\bf Proof of Theorem C.} By Theorem B and \rlemma{estfunct2}, the sequence $\{u^\Gw_{k,0}\}$ is nondecreasing and bounded from above by either $\Gl_1|x|^{-\gb_1}$ or $\Gl_2|x|^{-\gb_2}$. Therefore $\{u^\Gw_{k,0}\}$ increases to a function $u^\Gw_{\infty,0}$. By regularity theory, $u^\Gw_{\infty,0}$ is a solution of \eqref{A0} vanishing on $\prt \Gw \sms \{0\}$. Moreover, since $u^\Gw_{\infty,0} \geq u^\Gw_{k,0}$ for every $k>0$, $tr_{\prt \Gw}(u^\Gw_{\infty,0})=(\{0\},0)$. If $v \in \CU^\Gw_0$ then by using the argument as in the proof of \rth{gen+M}, we deduce that $v \geq u^\Gw_{\infty,0}$. Hence $u^\Gw_{\infty,0}=\min \CU^\Gw_0$. \qed \medskip

For any $\ell>0$ and any solution of \eqref{A0}, define
\bel{T} \Gw^\ell={\ell^{-1}}\Gw, \q T^1_\ell[u](x)=\ell^{\gb_1}u(\ell x), \q  T^2_\ell[u](x)=\ell^{\gb_2}u(\ell x) \q \forall x \in \Gw^\ell. \ee
\bprop{halfspace} Let $v \in C(\ovl \Gw \sms \{0\}) \cap C^2(\Gw)$ be a nonnegative solution of \eqref{A0} vanishing on $\prt \Gw \sms \{0\}$. \smallskip

\noindent{1)} Assume $H$ satisfies \eqref{multi}. For each $\ell$, put $v_\ell(x)=T^1_\ell[v](x)$. Then, up to a subsequence, $\{v_\ell\}$ converges in $C^1_{loc}(\ovl {\BBR_+^N} \sms \{0\})$, as $\ell \to 0$, to a solution of
\bel{half1} -\Gd u +  u^p |\nabla u|^q = 0 \text{ in } \BBR^N_+, \q u = 0 \text{ on } \prt \BBR^N_+\sms\{0\}. \ee
\noindent{2)} Assume $H$ satisfies \eqref{add}. For each $\ell$, put $v_\ell(x)=T^2_\ell[v](x)$. \smallskip

\noindent{2.i)} If $p=\frac{q}{2-q}$ then , up to a subsequence, $\{v_\ell\}$ converges in $C^1_{loc}(\ovl {\BBR_+^N} \sms \{0\})$, as $\ell \to 0$, to a solution of
\bel{half2} -\Gd u +   u^p + |\nabla u|^q = 0 \text{ in } \BBR^N_+, \q u = 0 \text{ on } \prt \BBR^N_+ \sms\{0\}. \ee

\noindent{2.ii)} If $p>\frac{q}{2-q}$ then , up to a subsequence, $\{v_\ell\}$ converges in $C^1_{loc}(\ovl {\BBR_+^N} \sms \{0\})$, as $\ell \to 0$, to a solution of
\bel{half3} -\Gd u +   u^p = 0 \text{ in } \BBR^N_+, \q u = 0 \text{ on } \prt \BBR^N_+\sms\{0\}. \ee

\noindent{2.iii)} If $p<\frac{q}{2-q}$ then , up to a subsequence, $\{v_\ell\}$ converges in $C^1_{loc}(\ovl {\BBR_+^N} \sms \{0\})$, as $\ell \to 0$, to a solution of
\bel{half4} -\Gd u +  |\nabla u|^q = 0 \text{ in } \BBR^N_+, \q u = 0 \text{ on } \prt \BBR^N_+\sms\{0\}. \ee
\es
\Proof We first notice that if $H$ satisfies either \eqref{multi} or \eqref{add} with $p=\frac{q}{2-q}$ then $v_\ell$ is a solution of \eqref{A0} in $\Gw^\ell$ which vanishes on $\prt \Gw^\ell \sms \{0\}$. If $H$ satisfies \eqref{add} with $p>\frac{q}{2-q}$ then $v_\ell$ satisfies
\bel{half3l} -\Gd v_\ell +   v_\ell^p + \ell^{\frac{p(2-q)-q}{p-1}}|\nabla v_\ell|^q= 0 \text{ in } \Gw^\ell, \q v_\ell = 0 \text{ on } \prt \Gw^\ell \sms \{0\}. \ee
If $H$ satisfies \eqref{add} with $p<\frac{q}{2-q}$ then $v_\ell$ satisfies
\bel{half4l} -\Gd v_\ell +  \ell^{\frac{q-(2-q)p}{q-1}}v_\ell^p + |\nabla v_\ell|^q= 0 \text{ in } \Gw^\ell, \q v_\ell = 0 \text{ on } \prt \Gw^\ell \sms \{0\}. \ee

Next, it follows from \rlemma{estfunct2} and \cite[Theorem 1]{Lib} that for every $R>1$ there exists a positive numbers $M=M(N,p,q,R)$ and $\gg=\gg(N,p,q) \in (0,1)$ such that
\bel{lie} \BA{ll} \sup\{\abs{v_\ell(x)}+\abs{\nabla v_\ell(x)}: x \in \Gg_{R^{-1},R} \cap \Gw^{\ell}\} \\[3mm]
   \phantom{qqqqq}
   + \sup\left\{\myfrac{\abs{\nabla v_\ell(x) - \nabla v_\ell (y)}}{\abs{x-y}^\gg}: x,y \in \Gg_{R^{-1},R} \cap \Gw^{\ell} \right\} \leq M \EA \ee
where $\Gg_{t_1,t_2}:=B_{t_2}(0)\sms B_{t_1}(0)$ with $0<t_1<t_2$. Notice that $M$ and $\gg$ are independent of $\ell \in (0,1)$ because the curvature of $\prt \Gw^{\ell}$ remains uniformly bounded when $0<\ell<1$. Thus there exists a sequence $\{\ell_n\}$ and a function $v^{\BBR^N_+} \in C^1(\ovl{\BBR^N_+} \sms \{0\})$ such that $\{v_{\ell_n}\}$ converges to $v^{\BBR^N_+}$ in $C^1_{loc}(\ovl{\BBR^N_+} \sms \{0\})$ which is a solution of
    \bel{halfspace} \left\{ \BA{lll} -\Gd v +  H(x,v,\nabla v) & = 0 \qq &\text{in } \BBR^N_+ \\
    \phantom{qerwerwhghqw,,,}
    v &= 0 &\text{in } \prt \BBR^N_+ \sms \{0\}\EA \right. \ee
Moreover,
    \bel{t2} \lim_{n \to \infty}(\sup\{|(v_{\ell_n}-v^{\BBR^N_+})(x)+|\nabla (v_{\ell_n}-v^{\BBR^N_+})(x)|: x \in \Gg_{R^{-1},R} \cap \Gw^{\ell_n}\})=0. \ee \qed

\bprop{halfspace-c} Under the assumptions of \rprop{halfspace}, if $v=u^\Gw_{\infty,0}$ then, up to a subsequence,  $\{v_\ell\}$ converge to $v^{\BBR^N_+}_{\infty}$  where $v^{\BBR^N_+}_{\infty}$  is a solution of
\bel{v-k}  \left\{ \BA{lll} (\ref{half1}) \text{ if } H \text{ satisfies } (\ref{multi}) \\
(\ref{half2}) \text{ if } H \text{ satisfies } (\ref{add}) \text{ with } p=\frac{q}{2-q} \\
(\ref{half3}) \text{ if } H \text{ satisfies } (\ref{add}) \text{ with } p>\frac{q}{2-q} \\
(\ref{half4}) \text{ if } H \text{ satisfies } (\ref{add}) \text{ with } p<\frac{q}{2-q}
\EA \right. \ee
with boundary trace $tr_{\prt \Gw}(v^{\BBR^N_+}_{\infty})=(\{0\},0)$.

Moreover, if $H$ satisfies \eqref{add}  then $ \{v_\ell\}$  converges to
$$
\left\{ \BA{ll} u^{\BBR^N_+}_{p,\infty,0} \text{ if } p > \frac{q}{2-q} \text{ where } u^{\BBR^N_+}_{p,\infty,0} \text{ is the solution of } (\ref{up-infty}) \text{ with } \Gw=\BBR^N_+,\\[3mm]
u^{\BBR^N_+}_{q,\infty,0} \text{ if } p < \frac{q}{2-q} \text{ where } u^{\BBR^N_+}_{q,\infty,0} \text{ is the solution of } (\ref{uq-infty})  \text{ with } \Gw=\BBR^N_+.
\EA \right. $$
\es
\Proof Since $v_\ell \geq u^{\Gw_\ell}_{k,0}$ for every $\ell>0$, $k>0$, $v^{\BBR^N_+} \geq u^{\BBR^N_+}_{k,0}$ for every $k>0$. Hence $tr_{\prt \Gw}(v^{\BBR^N_+})=(\{0\},0)$.

If $H$ satisfies \eqref{add} with $p \neq \frac{q}{2-q}$ then by uniqueness of strongly singular solution (see \cite{MV4} and \cite{NV}) either $\CU^{\BBR^N_+}_0=\{u^{\BBR^N_+}_{p,\infty,0}\}$ if $p>\frac{q}{2-q}$ or $\CU^{\BBR^N_+}_0=\{u^{\BBR^N_+}_{q,\infty,0}\}$ if $p<\frac{q}{2-q}$ . Therefore either $v^{\BBR^N_+}=u^{\BBR^N_+}_{p,\infty,0}$  if $p>\frac{q}{2-q}$ or $v^{\BBR^N_+}=u^{\BBR^N_+}_{q,\infty,0}$  if $p<\frac{q}{2-q}$ . \qed \medskip

We next study structure of the classes $\CE_i$ ($i=1,2,3,4$) ($\CE_i$ is defined in \eqref{PHi}).

\bth{exis-uniq-SS}  i)  Assume either $H$ satisfies \eqref{multi} with $N(p+q-1) \geq p+1$ or $H$ satisfies \eqref{add} with $m_{p,q} \geq p_c$. Then $\CE_i = \emptyset$  where $i=1$ if $H$ satisfies \eqref{multi}, $i=2$ if $H$ satisfies \eqref{add} with $p=\frac{q}{2-q}$, $i=3$ if $H$ satisfies \eqref{add} with $p>\frac{q}{2-q}$, $i=4$ if $H$ satisfies \eqref{add} with $p<\frac{q}{2-q}$.

\noindent ii) Assume either $H$ satisfies \eqref{multi} with $0<N(p+q-1)<p+1$ or $H$ satisfies \eqref{add} with $m_{p,q}<p_c$. Then $\CE_i \neq \emptyset$ with $i \in \ovl {1,4}$ determined  in the statement i).
\es
\Proof Notice that when  $H$ satisfies \eqref{add} with $p \neq \frac{q}{2-q}$, the statement i) and ii) have been proved in \cite{MV4} and \cite{NV}. Moreover, if $m_{p,q}=p_c$ then $\CE_i=\{\gw_i^*\}$ with $i=3$ if $p>\frac{q}{2-q}$, $i=4$ if $p<\frac{q}{2-q}$.  So we are left with the case when $H$ satisfies either \eqref{multi} or $H$ satisfies \eqref{add} with $p=\frac{q}{2-q}$ and we only give the proof for the case $H$ satisfies \eqref{multi}.

\noindent {\it Proof of statement i).} Denote by $\vgf_1$ the first eigenfunction of $-\Gd'$ in $W^{1,2}_0(S^{N-1}_+)$, normalized such that $\max_{S^{N-1}_+}\vgf_1=1$, with corresponding eigenvalue $\gl_1=N-1$. Multiplying \eqref{PH1} by $\vgf_1$ and integrating over $S_+^{N-1}$, we get
	$$ \BA{ll}
	\left[N-1-\gb_1(\gb_1+2-N)\right]\myint{S_+^{N-1}}{}\gw\, \vgf_1 dx \\ \phantom{-----}
+ \myint{S_+^{N-1}}{}\gw^p(\gb_1^2\, \gw^2+\abs{\nabla' \gw}^2)^{\frac{q}{2}}\vgf_1 dx = 0. \EA $$
{\it Therefore  if $N-1 \geq \gb_1(\gb_1+2-N)$, namely $N(p+q-1)\geq p+1$, then there exists no positive solution of \eqref{PH1}}. \medskip

\noindent{\it Proof of statement ii).} The proof is based on construction a subsolution and a supersolution to \eqref{PH1}. By a computation, we can prove that $\unl \gw:=\gth_1 \vgf_1^{\gth_2}$ is a positive subsolution of \eqref{PH1} with $\gth_1>0$ small and $ 1< \gth_2< \frac{\gb_1\,(\gb_1+2-N)}{N-1}$. Next, it is easy to see that $\ovl \gw=\gth_4$, with $\gth_4>0$ large enough, is a supersolution of \eqref{PH1} and $\ovl \gw>\unl \gw$ in $\ovl S^{N-1}_+$. Therefore by \cite{KaKr} there exists a solution $\gw^*_1 \in W^{2,m}(S^{N-1}_+)$ (for any $m>N$) to \eqref{PH1} such that $0< \unl \gw \leq \gw^*_1 \leq \ovl \gw$ in $S_+^{N-1}$. By regularity theory,  $\gw^*_1 \in C^{\infty}(\overline{S^{N-1}_+})$.  \medskip \qed

The structure of $\CE_i$ ($i=1,2 $) is analyzed in the following theorem.

\bth{unique-sphere} \noindent  i) If $H$ satisfies \eqref{multi} with $N(p+q-1)<p+1$ and $p \geq 1$ then $\CE_1=\{\gw_1^*\}$.

\noindent ii) If $H$ satisfies \eqref{add} with $m_{p,q}<p_c$ and $p=\frac{q}{2-q}$ then $\CE_2=\{\gw_2^*\}$.
\es
\Proof We give below only the proof of statement i) because the statement ii) can be treated in a similar way. Suppose that $\gw_1$ and $\gw_2$ are two positive different solutions of \eqref{PH1} and by Hopf lemma $\nabla'\gw_i$ ($i=1,2$) does not vanish on $S^{N-1}_+$.  Up to exchanging the role of $\gw_1$ and $\gw_2$, we may assume $ \max_{S^{N-1}_+}\gw_2 \geq \max_{S^{N-1}_+}\gw_1$ and
	$$ \gt_0:=\inf\{\gt>1:\gt\gw_1 > \gw_2 \text{ in } S^{N-1}_+\} >  1. $$
Set $\gw_{1,\gt_0}:=\gt_0 \gw_1$, then $\gw_{1,\gt_0}$ is a positive supersolution to problem \eqref{PH1}.Put $\tl \gw=\gw_{1,\tau_0}-\gw_2 \geq 0$. If there exists $\gs_0 \in S^{N-1}_+$ such that $\gw_{1,\gt_0}(\gs_0)=\gw_2(\gs_0)>0$ and $\nabla ' \gw_{1,\gt_0}(\gs_0)=\nabla ' \gw_2(\gs_0)$ then $\tl \gw(\gs_0)=0$ and $\nabla' \tl \gw(\gs_0)=0$. This contradicts the strong maximum principle (see \cite{GT}). If $\gw_{1,\gt_0}>\gw_2$ in $S^{N-1}_+$ and  there exists $\gs_0 \in \prt S^{N-1}_+$ such that $\frac{\prt \gw_{1,\gt_0}}{\prt \gn}(\gs_0)=\frac{\prt \gw_2}{\prt \gn}(\gs_0)$ then $\tl \gw > 0$ and $\frac{\prt \tl \gw}{\prt \gn}(\gs_0)=0$. This contradict the Hopf lemma (see \cite{GT}). \qed \medskip

When $\BBR^N_+$ is replaced by a general $C^2$ bounded domain $\Gw$, the role of $\gw^*_i$ is crucial for describing the strong singularities. In that case we assume that $0\in\prt\Gw$ and the tangent plane to $\prt\Gw$ at $0$ is $\prt\BBR^{N-1}_+:=\{(x',0):x'\in \BBR^{N-1}\}$, with normal inward unit vector ${\bf e}_N$.

Let $u \in C(\ovl {\BBR^N_+}\sms\{0\})$ be a solution of \eqref{P3}. When $H$ satisfies \eqref{multi}, $T^1_\ell[u]$ is  a solution of \eqref{P3} and we say that $u$ is {\em self-similar} if $T^1_\ell[u] = u$ for every $\ell$. When $H$ satisfies \eqref{add} with $p=\frac{q}{2-q}$,  $T^2_\ell[u]$ is  a solution of \eqref{P3} and we say that $u$ is {\em self-similar} if $T^2_\ell[u] = u$ for every $\ell>0$.

\bprop{sing}

\noindent i) If $H$ satisfies \eqref{multi} with $N(p+q-1)<p+1$ and $p \geq 1$ then
\bel{uniq7}
\lim_{\tiny\BA{c}\Gw \ni x\to 0\\
\frac{x}{|x|}=\gs\in S^{N-1}_+
\EA}|x|^{\gb_1}u^\Gw_{\infty,0}(x)=\gw_1^*(\gs),
 \ee
 locally uniformly on $S^{N-1}_+$. \smallskip

\noindent ii) If $H$ satisfies \eqref{add} with $m_{p,q}<p_c$ then
\bel{uniq7'}
\lim_{\tiny\BA{c}\Gw \ni x\to 0\\
\frac{x}{|x|}=\gs\in S^{N-1}_+
\EA}|x|^{\gb_2}u^\Gw_{\infty,0}(x)=\gw_i^*(\gs),
 \ee
locally uniformly on $S^{N-1}_+$ where $i=2$ if $p=\frac{q}{2-q}$, $i=3$ if $p>\frac{q}{2-q}$, $i=4$ if $p<\frac{q}{2-q}$.\smallskip
\es
\Proof {\bf Case 1: $H$ satisfies \eqref{multi}.} Since the proof is close to the one of \cite[Proposition 3.22]{NV}, we present the main ideas.

We first  note that $T^1_\ell[u^{\BBR^N_+}_{\infty,0}]=u^{\BBR^N_+}_{\infty,0}$ for every $\ell>0$. Hence $u^{\BBR^N_+}_{\infty,0}$ is self-similar and satisfies \eqref{uniq7} with $\Gw$ replaced by $\BBR^N_+$.

Next,  let $B$ and $B'$ are two open balls tangent to $\prt \Gw$ at $0$ such that $B \sbs \Gw \sbs G:=(B')^c$. Then
 \bel{uniq12*} u_{\infty,0}^{B^{\ell'}}\leq u_{\infty,0}^{B^{\ell}} \leq u_{\infty,0}^{\BBR^N_+}\leq  u_{\infty,0}^{G^\ell}
 \leq  u_{\infty,0}^{G^{\ell ''}} \qq\forall \,0<\ell\leq\ell',\ell''\leq 1.
\ee
Notice that $u_{\infty,0}^{B^{\ell}}\uparrow \underline u_{\infty,0}^{\BBR^N_+}$ and $u_{\infty,0}^{G^\ell}\downarrow \overline u_{\infty,0}^{\BBR^N_+}$ when $\ell\to 0$ where $ \underline u_{\infty,0}^{\BBR^N_+}$ and $\overline u_{\infty,0}^{\BBR^N_+}$ are positive solutions of \eqref{A0} in $\BBR^N_+$, continuous in $\overline{\BBR^N_+}\setminus\{0\}$ and vanishing on $\prt\BBR^N_+\setminus\{0\}$. By letting $\ell \to 0$ in \eqref{uniq12}, we obtain
  \bel{uniq13a}
 u_{\infty,0}^{B^{\ell}} \leq \underline u_{\infty,0}^{\BBR^N_+}\leq u_{\infty,0}^{\BBR^N_+}\leq \overline u_{\infty,0}^{\BBR^N_+}\leq  u_{\infty,0}^{G^\ell}\qq\forall\,0<\ell\leq 1.
\ee
Furthermore there also holds for $\ell,\ell'>0$,
  \bel{uniq14}
T^1_{\ell'\ell}[u_{\infty,0}^{B}]=T^1_{\ell'}[T^1_{\ell}[u_{\infty,0}^{B}]]=u_{\infty,0}^{B^{\ell\ell'}} \text{ and }
T^1_{\ell'\ell}[u_{\infty,0}^{G}]=T^1_{\ell'}[T^1_{\ell}[u_{\infty,0}^{G}]]=u_{\infty,0}^{G^{\ell\ell'}}.
\ee
Letting $\ell\to 0$ in \eqref{uniq14} yields
  \bel{uniq15}
\underline u_{\infty,0}^{\BBR^N_+}=T^1_{\ell'}[\underline u_{\infty,0}^{\BBR^N_+}] \text{ and }
\overline u_{\infty,0}^{\BBR^N_+}=T^1_{\ell'}[\overline u_{\infty,0}^{\BBR^N_+}].
\ee
Thus $\underline u_{\infty,0}^{\BBR^N_+}$ and $\overline u_{\infty,0}^{\BBR^N_+}$ are self-similar solutions of \eqref{A0} in $\BBR^N_+$ vanishing on $\prt \BBR_+^N\sms\{0\}$ and continuous in $\overline{\BBR^N_+}\setminus\{0\}$. Therefore they coincide with $u_{\infty,0}^{\BBR^N_+}$.

Finally, since
 \bel{uniq16} u_{\infty,0}^{B^{\ell}}\leq T^1_{\ell}[u_{\infty,0}^{\Gw}]
 \leq  u_{\infty,0}^{G^\ell}\qq\forall \,0<\ell\leq 1,
\ee
by letting $\ell \to 0$ we obtain \eqref{uniq7}. \medskip

\noindent {\bf Case 2: $H$ satisfies \eqref{add} with $p=\frac{q}{2-q}$.} The proof is similar to the one in case 1.

\noindent {\bf Case 3: $H$ satisfies \eqref{add} with $p>\frac{q}{2-q}$.} For any $k>0$ and $\ell>0$, $T^2_\ell[u^\Gw_{k,0}]$ is a solution of  \eqref{half3l} with boundary trace $k\gd_0$. Denote by $u^{\Gw^\ell}_{p,k,0}$ the solution of
\bel{upk} -\Gd u +   u^p = 0 \text{ in } \Gw^\ell, \q u=k\gd_0 \text{ on } \prt \Gw^\ell. \ee
Since $0<\ell<1$ and $p>\frac{q}{2-q}$, by comparison principle, we get
$$u^{\Gw^\ell}_{\ell^{\gb_2+1-N} k,0} \leq T^2_\ell[u^\Gw_{k,0}] \leq u^{\Gw^\ell}_{p,\ell^{\gb_2+1-N}k,0}.$$
in $\Gw^\ell$. By letting $k \to \infty$, we obtain
$$u^{\Gw^\ell}_{\infty,0} \leq T^2_\ell[u^\Gw_{\infty,0}] \leq u^{\Gw^\ell}_{p,\infty,0}$$
in $\Gw^\ell$ where  $u^{\Gw^\ell}_{p,\infty,0}$ is the unique problem of
\bel{up-inf} -\Gd u +   u^p = 0 \text{ in } \Gw^\ell, \q tr_{\prt \Gw}(u)=(\{0\},0) \text{ on } \prt \Gw^\ell. \ee
By \rprop{halfspace-c}, letting $\ell \to 0$ we deduce that
$$ \lim_{\ell \to 0}\ell^{\gb_2}u^\Gw_{\infty,0}(\ell x)=u^{\BBR^N_+}_{p,\infty,0}(x) $$
where $u^{\BBR^N_+}_{p,\infty,0}$ is the unique solution of \eqref{half3} with strong singularity at $0$. Hence, it follows from \cite{MV4} that $u^\Gw_{\infty,0}$ satisfies \eqref{uniq7'} with $i=3$. \medskip

\noindent {\bf Case 4: $H$ satisfies \eqref{add} with $p<\frac{q}{2-q}$.} By similar argument in case 3 and results in \cite{NV}, we derive \eqref{uniq7'} with $i=4$. \qed \medskip

We next construct the maximal strongly singular solution at $0$.

\bprop{max}  i) Assume either $H$ satisfies \eqref{multi} with $0<N(p+q-1)<p+1$ then there exists a maximal element $U_{\infty,0}^{\Gw}$ of $\CU^\Gw_0$. In addition, if $p \geq 1$ then
\bel{max1}
\lim_{\tiny\BA{c}\Gw \ni x\to 0\\
\frac{x}{|x|}=\gs\in S^{N-1}_+
\EA}|x|^{\gb_1}U_{\infty,0}^{\Gw}(x)=\gw^*_1(\gs),
 \ee
 locally uniformly on $S^{N-1}_+$. \medskip

\noindent ii) If $H$ satisfies \eqref{multi} with $m_{p,q}<p_c$ then there exists a maximal element $U_{\infty,0}^{\Gw}$ of $\CU^\Gw_0$ and
\bel{max1-add}
\lim_{\tiny\BA{c}\Gw \ni x\to 0\\
\frac{x}{|x|}=\gs\in S^{N-1}_+
\EA}|x|^{\gb_2}U_{\infty,0}^{\Gw}(x)=\gw^*_i(\gs),
 \ee
 locally uniformly on $S^{N-1}_+$ where $i=2$ if $p=\frac{q}{2-q}$, $i=3$ if $p>\frac{q}{2-q}$, $i=4$ if $p<\frac{q}{2-q}$.
\es
\Proof {\bf Case 1:  $H$ satisfies \eqref{multi}.} \smallskip

\noindent {\it Step 1: Construction maximal solution}. Let $u$ is a positive solution of \eqref{A0} which vanishes on $\prt \Gw \sms \{0\}$. Since $0<N(p+q-1)<p+1$, there exists a radial solution of \eqref{A0} in $\BBR^N\setminus\{0\}$ of the form
\bel{max2}
U^S_1(x)=\Gl^\Gw_3\,|x|^{-\gb_1}\q\text{with }\;\Gl^S_1=\left(\myfrac{\gb_1+2-N}{\gb_1^{q-1}}\right)^{\frac{1}{p+q-1}}.
\ee
Therefore,
 $U_1^{\Gw,*}(x)=\Gl^*_1|x|^{-\gb_1}$ with $\Gl^*_1= \max\{\Gl^S_1,\Gl_1\}$ is a supersolution of \eqref{A0} in $\BBR^N\setminus\{0\}$ and dominates $u$ in $\Gw$. For $0<\ge<\max\{|z|:z\in\Gw\}$, we construct a decreasing smooth sequence $\{\psi_{\ge,n}\}$ on $(\prt \Gw \sms B_\ge(0)) \cup (\Gw \cap \prt B_\ge(0))$ as follows
$$ 0\leq \psi_{\ge,n} \leq \Gl^*_1\ge^{-\gb_1}, \quad  \psi_{\ge,n}(x)=\Gl^*_1\ge^{-\gb_1}\quad \text{ if } x \in \Gw \cap \prt B_\ge(0) $$
$$ \psi_{\ge,n}(x)=0 \quad \text{if } x \in \prt \Gw \sms B_\ge(0) \text{ and } \dist(x,\prt B_\ge(0))>\frac{1}{n}. $$
Let $u^\Gw_{\ge,n}$ the solution of
 \bel{u aprox}
\left\{\BA {lll}
-\Gd u+ u^p|\nabla u|^q&=0\qq&\text {in }\Gw\setminus B_\ge(0)\\\phantom{-\Gd +|\nabla u_\ge|^q,,}
u&=\psi_{\ge,n} \qq&\text {on }(\prt\Gw\setminus B_\ge(0)) \cup (\Gw \cap \prt B_\ge(0))
\EA\right.\ee
By the comparison principle, $ u \leq u_{\ge,n} \leq U^{\Gw,*}_1 $ in $\Gw\sms B_\ge(0)$. Owing to \rcor{var-stab}, $\{u_{\ge,n}\}$ converges to the solution $u^\Gw_\ge$ of
  \bel{max4} \left\{\BA {lll}
-\Gd u_\ge+u_\ge^p|\nabla u_\ge|^q&=0\qq&\text {in }\Gw\setminus B_\ge(0) \\ \phantom{----- ,,,,} 
u_\ge&=0\qq&\text {on }\prt\Gw\setminus B_\ge(0)\\\phantom{------,,}
u_\ge&=\Gl^*_1\ge^{-\gb_1}\qq&\text {on }\Gw\cap \prt B_\ge(0).
\EA\right.\ee
Consequently, $ u \leq u^\Gw_{\ge} \leq U^{\Gw,*}_1$. If  $\ge'<\ge$, for $n$ large enough, $u^\Gw_{\ge',n} \leq u^\Gw_{\ge,n}$, therefore
  \bel{max5}
  u\leq u^\Gw_{\ge'}\leq u^\Gw_\ge\leq U^{\Gw,*}_1(x)\qq\text{in }\;\Gw.
  \ee
Letting $\ge$ to zero, $\{u^\Gw_\ge\}$ decreases and converges to some $U_{\infty,0}^{\Gw}$ which vanishes on  $\prt\Gw\setminus \{0\}$.   By regularity theory, the convergence occurs in $C^1_{loc}(\overline\Gw\setminus\{0\})$, $U_{\infty,0}^{\Gw}\in \CU^\Gw_0$. Moreover, there holds
    \bel{max6}
 u_{\infty,0}^{\Gw} \leq  u\leq U_{\infty,0}^{\Gw}\leq U^{\Gw,*}_1(x).
  \ee
Therefore $ U_{\infty,0}^{\Gw}$ is the maximal element of $\CU^\Gw_0$.\

Notice that when $H$ satisfies \eqref{add} with $p>\frac{q}{2-q}$, there is no radial solution of \eqref{A0} in $\BBR^N \sms \{0\}$. We can instead employ a radial supersolution of the form
\bel{max2-add}
U^S_3(x)=\Gl^\Gw_3\,|x|^{-\gb_2}\q\text{with }\;\Gl^S_3={\gb_2(\gb_2+2-N)}^{\frac{1}{p-1}}
\ee
and then we proceed as above to construct the maximal solution.\medskip

  \noindent{\it Step 2: Proof of \eqref{max1}.} Assume $H$ satisfies \eqref{multi} with $p \geq 1$. We first take into account that $U_{\infty,0}^{\BBR^N_+} \equiv u_{\infty,0}^{\BBR^N_+}$. Indeed, the assertion follows from the fact that
    \bel{max7} T^1_\ell[U^{\Gw,*}_1]|_{_{|x|=\ge}}=U^{\Gw,*}_1|_{_{|x|=\ge}}\qq\forall\,\ell>0,
    \ee
and the \rth{unique-sphere}.

Next,  let $B$ and $B'$ are two open balls tangent to $\prt \Gw$ at $0$ such that $B \sbs \Gw \sbs G:=(B')^c$. Note that $T^1_\ell[u^\Gth_\ge]=u^{\Gth^\ell}_{\frac{\ge}{\ell}}$ for any $\ell,\ge>0$ and any domain $\Gth$ (with $0\in \prt\Gth$) where $u^\Gth_\ge$ is the solution of \eqref{max4} in $\Gth \setminus B_\ge(0)$. By taking $\Gth=B$ and $\Gth=G$ successively and by letting $\ge \to 0$ we deduce that
\bel{max9}
 T^1_\ell[U_{\infty,0}^{B}]=U_{\infty,0}^{B^{\ell}} \text{ and }\;T^1_\ell[U_{\infty,0}^{G}]=U_{\infty,0}^{G^\ell}.
\ee
By comparison,
  \bel{uniq12} U_{\infty,0}^{B^{\ell'}}\leq U_{\infty,0}^{B^{\ell}} \leq U_{\infty,0}^{\BBR^N_+}\leq  U_{\infty,0}^{G^\ell}
 \leq  U_{\infty,0}^{G^{\ell''}}\qq\forall \,0<\ell\leq\ell',\ell''\leq 1
\ee
and
  \bel{uniq13} U_{\infty,0}^{B^{\ell'}}\leq U_{\infty,0}^{B^{\ell}} \leq T^1_\ell[U_{\infty,0}^{\Gw}] \leq  U_{\infty,0}^{G^\ell}
 \leq  U_{\infty,0}^{G^{\ell''}}\qq\forall \,0<\ell\leq\ell',\ell''\leq 1.
\ee
Hence $U_{\infty,0}^{B^{\ell}}\uparrow \underline U_{\infty,0}^{\BBR^N_+}\leq U_{\infty,0}^{\BBR^N_+}$ and
$U_{\infty,0}^{G^\ell} \downarrow \overline U_{\infty,0}^{\BBR^N_+}\geq U_{\infty,0}^{\BBR^N_+}$ as $\ell \to 0$ where  $\underline U_{\infty,0}^{\BBR^N_+}$ and $ \overline U_{\infty,0}^{\BBR^N_+}$ are positive solutions of \eqref{A0} in $\BBR^N$ which vanish on $\prt\BBR^N_+\setminus\{0\}$ and endow the same scaling invariance under $T^1_\ell$. Therefore they coincide with $u_{\infty,0}^{\BBR^N_+}$. Letting $\ell \to 0$ in \eqref{uniq13} implies \eqref{max1}.  \medskip

\noindent {\bf Case 2: $H$ satisfies \eqref{add} with $p=\frac{q}{2-q}$.} The proof is smilar to the one in case 1. \medskip

\noindent {\bf Case 3: $H$ satisfies \eqref{add} with $p>\frac{q}{2-q}$.} Since $u^{\Gw^\ell}_{\frac{\ge}{\ell}} \leq T^2_\ell[u^\Gw_\ell]$, by letting $\ge \to 0$ we obtain $U^{\Gw^\ell}_{\infty,0} \leq T^2_\ell[U^\Gw_{\infty,0}]$. It follows that
$$ u^{\Gw^\ell}_{\infty,0} \leq  U^{\Gw^\ell}_{\infty,0}  \leq  T^2_\ell[U^\Gw_{\infty,0}] \leq  T^2_\ell[u^\Gw_{p,\infty,0}] =  u^{\Gw^\ell}_{p,\infty,0} $$
where $u^{\Gw^\ell}_{p,\infty,0}$ is the solution of \eqref{up-inf}. Due to \rprop{halfspace-c} and the uniqueness, we deduce
	$$ \lim_{\ell \to 0}T^2_\ell[U^\Gw_{\infty,0}]=u^{\BBR^N_+}_{p,\infty,0},$$
from which \eqref{max1-add} follows straightforward. \medskip

\noindent {\bf Case 4: $H$ satisfies \eqref{add} with $p<\frac{q}{2-q}$.} The proof is similar to the one in case 3. \qed \medskip

 \rprop{sing} and \rprop{max} show that the minimal solution $u^\Gw_{\infty,0}$ and the maximal solution $U^\Gw_{\infty,0}$ have the same asymptotic behavior near $0$, which allows us to prove the following result.

\bth {UNI} Assume either $H$ satisfies \eqref{multi} with $N(p+q-1)<p+1$ and $p \geq 1$ or $H$ satisfies \eqref{add} with $m_{p,q}<p_c$. Then $U^\Gw_{\infty,0}=u^\Gw_{\infty,0}$.
  \es
\Proof We follow the method used in \cite[Sec 4]{GV}. \medskip

\noindent {\bf Case 1: $H$ satisfies \eqref{multi} with $p\geq 1$.} \smallskip

We represent $\prt\Gw$ near $0$ as the graph of a $C^2$ function $\gf$
defined in $\BBR^{N-1}\cap B_R$ and such that $\gf(0)=0$, $\nabla_{N-1}\gf(0)=0$ and
$$\prt\Gw\cap B_R=\{x=(x',x_N):x'\in \BBR^{N-1}\cap B_R,x_N=\gf(x')\}.$$
We introduce the new variable $y=\Gf(x)$ with $y'=x'$ and $y_N=x_N-\gf(x')$, with corresponding spherical coordinates in $\BBR^N$, $(r,\gs)=(|y|,\frac{y}{|y|})$.

Let $u$ is a positive solution of \eqref{A0} in $\Gw$ vanishing on $\prt\Gw \sms \{0\}$. We set $u(x)=r^{-\gb_1}v(t,\gs)$ with  $t=-\ln r\geq 0$, then a technical computation shows that $v$ satisfies with ${\bf n}=\frac{y}{|y|}$
\begin{equation}\label{uni2}\BA {l}
\left(1+\ge^1_1\right)v_{tt}+\left(2\gb_1+2-N+\ge^1_2\right)v_{t}
+(\gb_1\left(\gb_1+2-N)+\ge^1_3\right)v+\Gd'v\\[3mm]
\phantom{--}
+\langle\nabla'v,\overrightarrow {\ge^1_4}\rangle+\langle\nabla'v_t,\overrightarrow {\ge^1_5}\rangle+
\langle\nabla'\langle \nabla' v,{\bf e}_N\rangle,\overrightarrow {\ge^1_6}\rangle\\[3mm]\phantom{--}
-v^p\abs{(-\gb_1\,v+v_t){\bf n}+\nabla' v+\langle(-\gb_1\,v+v_t){\bf n}+\nabla' v,{\bf e}_N\rangle\overrightarrow\ge^1_7}^q=0,
\EA\end{equation}
on $Q_R:=[-\ln R, \infty)\ti S^{N-1}_{+}$ and vanishes on $[-\ln R,\infty)\ti \prt S^{N-1}_{+}$, where
Furthermore the $\ge^1_j$ are uniformly continuous functions of $t$ and $\gs\in S^{N-1}$ for $j=1,...,7$, $C^1$ for $j=1,5,6,7$ and satisfy the following decay estimates $|\ge^1_j(t,.)|\leq c_{12}e^{-t}$  for $j=1,...,7$ and $|\ge^1_{j\,t}(t,.)|+|\nabla'\ge^1_j|\leq c_{12}e^{-t}$  for $j=1,5,6,7$.
By \cite[Theorem 4.7]{GV}, there exist a constant $c_{13}>0$ and $T>\ln R$ such that
\begin{equation}\label{uni4}\BA {l}
\norm {v(t,.)}_{C^{2,\gg}(\overline {S^{N-1}_{+}})}+\norm {v_t(t,.)}_{C^{1,\gg}(\overline {S^{N-1}_{+}})}
+\norm {v_{tt}(t,.)}_{C^{0,\gg}(\overline {S^{N-1}_{+}})}\leq c_{13}
\EA\end{equation}
for any $\gamma\in (0,1)$ and $t \geq T+1$. Moreover $ \lim_{t \to \infty}^{}\int_{S^{N-1}_+}^{}(v_t^2 + v_{tt}^2 + |\nabla ' v_t|^2)d\gs = 0$. Consequently, the $\gw$-limit set of $v$
$$ \Gg^+(v)=\cap_{\tau \geq 0}\ovl{\cup_{t \geq \tau}v(t,.)}^{C^2(S^{N-1}_+)} $$
is a non-empty, connected and compact subset of the set of $\CE_1$. By the uniqueness of \eqref{PH1}, $\Gg^+(v)=\CE_1=\{\gw^*_1\}$. Hence $\lim_{t\to \infty}v(t,.)=\gw^*_1$ in $C^2(\overline {S^{N-1}_{+}})$.\smallskip

By taking $u=u^\Gw_{\infty,0}$ and $u=U^\Gw_{\infty,0}$ we obtain
\bel{uni5}
\lim_{\Gw \ni x\to 0}\myfrac{u^\Gw_{\infty,0}(x)}{U^\Gw_{\infty,0}(x)}=1.
\ee
For any $\vge>0$, by the comparison principle, $(1+\vge)u_{\infty,0}^\Gw \geq U_{\infty,0}^\Gw$ in $\Gw \sms B_\vge(0)$. Letting $\vge \to 0$ yields $u_{\infty,0}^\Gw \geq U_{\infty,0}^\Gw$ in $\Gw$ and thus $u_{\infty,0}^\Gw = U_{\infty,0}^\Gw$ in $\Gw$. \medskip

\noindent{\bf Case 2: $H$ satisfies \eqref{add} with $p=\frac{q}{2-q}$.} The assertion is obtained by a similar argument. \medskip

\noindent{\bf Case 3: $H$ satisfies \eqref{add} with $p>\frac{q}{2-q}$.}
In this case, we use the transformation $t=-\ln r$ for $t\geq 0$ and $ \tilde u(r,\gs)=r^{-\gb_2}v(t,\gs)$ and obtain the following equation instead of \eqref{uni2}
\bel{uni2-add}\BA{lll}
\left(1+\ge^3_1\right)v_{tt}+\left(2\gb_2+2-N+\ge^3_2\right)v_{t}
+\left(\gb_2(\gb_2+2-N)+\ge^3_3\right)v+\Gd'v\\[3mm]
+\langle\nabla'v,\overrightarrow {\ge^3_4}\rangle+\langle\nabla'v_t,\overrightarrow {\ge^3_5}\rangle+
\langle\nabla'\langle \nabla' v,{\bf e}_N\rangle,\overrightarrow {\ge^3_6}\rangle -v^p \\[3mm]
- e^{-\frac{p(2-q)-q}{p-1}t}\abs{(-\gb_1\,v+v_t){\bf n}+\nabla' v+\langle(-\gb_1\,v+v_t){\bf n}+\nabla' v,{\bf e}_N\rangle\overrightarrow\ge^3_7}^q=0
\EA\ee
where $\ge^3_j$ has the same properties as $\ge^1_j$ ($j=\ovl{1,7}$). Notice that $$\lim_{t \to \infty}e^{-\frac{p(2-q)-q}{p-1}t}=0$$ since $p>\frac{q}{2-q}$. By proceeding as in the Case 1, we deduce that $u_{\infty,0}^\Gw = U_{\infty,0}^\Gw$ in $\Gw$. \medskip

\noindent{\bf Case 4: $H$ satisfies \eqref{add} with $p<\frac{q}{2-q}$.} Using a similar argument as in Case 3, we obtain $u_{\infty,0}^\Gw = U_{\infty,0}^\Gw$ in $\Gw$. \qed \medskip

\noindent{\bf Proof of Theorem D.} The proof follows by combining \rth{exis-uniq-SS}, \rth{unique-sphere}, \rprop{sing} and \rth{UNI}. \qed \medskip

\noindent{\bf Proof of Theorem E.} \smallskip

\noindent {\bf Case 1: $H(.,u,\nabla u) \in L^1_\rho(\Gw)$.} It follows from \rprop{tra-reg} that $\CR(u)=\prt \Gw$ and hence the boundary trace of $u$ is a bounded Radon measure on $\prt \Gw$. Since $u=0$ on $\prt \Gw \sms \{0\}$, $\gm=k\gd_0$ for some $k \geq 0$. If $k=0$, then $u \equiv 0$. If $k>0$ then by \rth{unique} $u=u^\Gw_{k,0}$ and \eqref{com} follows from \rprop{compare}. \medskip

\noindent{\bf Case 2: $H(.,u,\nabla u) \notin L^1_\rho(\Gw)$.} By \rth{gen+M}, $tr_{\prt \Gw}(u)=(\{0\},0)$. From \rth{UNI}, $u=u^\Gw_{\infty,0}$. \qed

\section{Removability}
In this section we deal with removable singularities in the case that $H$ is critical or supercritical.
\bprop{nostrong} Assume either $H$ satisfies \eqref{multi} with $N(p+q-1)> p+1$
or H satisfies \eqref{add} with $m_{p,q} > p_c$. If  $u\in C(\overline\Gw\setminus\{0\})\cap C^2(\Gw)$ is a nonnegative solution of \eqref{A0} vanishing on $\prt\Gw\setminus\{0\}$ then $u$ cannot be strongly singular solution. \es
\Proof  We consider a sequence of functions $\gz_n\in C^\infty(\BBR^N)$ such that $\gz_n(x)=0$ if $|x|\leq \frac{1}{n}$, $\gz_n(x)=1$ if $|x|\geq \frac{2}{n}$, $0\leq \gz_n\leq 1$ and $|\nabla\gz_n|\leq c_{13}n$, $|\Gd\gz_n|\leq c_{13}n^2$ where $c_{13}$ is independent of $n$. We take $\gx\gz_n$ as a test function (where $\gx$ is the solution to \eqref{eta}) and we obtain
\bel{est3}\BA {l}
\myint{\Gw}{}( u + H(x,u,\nabla u)\gx)\gz_n \, dx=\myint{\Gw}{}u\left(\gx\Gd \gz_n+2\nabla\gx.\nabla\gz_n \right)dx
=J_1+J_2.
\EA\ee
Set $\CO_n=\Gw\cap \{x:\frac{1}{n}<|x|\leq \frac{2}{n}\}$, then $|\CO_n|\leq c_{14}(N)n^{-N}$. On one hand
$$  J_1 \leq c_{15}\Gl_i\,\myint{\CO_n}{}n^{\gb_i+2}\gx dx\leq c_{16}n^{\gb_i+1-N}
$$
since $\gx(x)\leq c_3\rho(x)$ where
\bel{i} i=\left\{ \BA{lll} 1 \q \text{if } H \text{ satisfies } (\ref{multi}), \\ 2 \q \text{if } H \text{ satisfies } (\ref{add}). \EA \right. \ee
On the other hand,
\bel{est3'} J_2 \leq c_{17}\Gl_i\,\myint{\CO_n}{}n^{\gb_i+1}|\nabla\gx| dx\leq c_{18}n^{\gb_i+1-N} \ee
where $i$ is given by \eqref{i}. By combining \eqref{est3}-\eqref{est3'} and then by letting $n \to \infty$ we obtain
\bel{est4}\BA {l}
\myint{\Gw}{}\left(u+ H(x,u,\nabla u)\gx\right) dx < \infty. \EA\ee
Hence $u$ is a moderate solution of \eqref{A0}. Therefore the boundary trace of $u$ is a bounded measure. Since $u=0$ on $\prt \Gw \sms \{0\}$, the boundary trace of $u$ is $k\gd_0$ for some $k \geq 0$. \qed

\bcor{supercritical} Assume either $H$ satisfies \eqref{multi} with $N(p+q-1) > p+1$
or H satisfies \eqref{add} with $m_{p,q} > p_c$. If  $u\in C(\overline\Gw\setminus\{0\})\cap C^2(\Gw)$ is a nonnegative solution of \eqref{A0} vanishing on $\prt\Gw\setminus\{0\}$ then $u \equiv 0$. \es
\Proof We deduce from the assumption, \eqref{est3}-\eqref{est3'} that
$$ \myint{\Gw}{}(u+H(x,u,\nabla u)\gx) dx = 0, $$
which implies $u \equiv 0$. \qed

\bth{critical} Assume either $H$ satisfies \eqref{multi} with $N(p+q-1) = p+1$
or H satisfies \eqref{add} with $m_{p,q} = p_c$. If $u\in C(\overline\Gw\setminus\{0\})\cap C^2(\Gw)$ is a nonnegative solution of \eqref{A0} vanishing on $\prt\Gw\setminus\{0\}$ then $u\equiv 0$.
\es
\Proof By \rprop{nostrong}, $u$ admits a boundary trace  $k\gd_0$ with $k \geq 0$.

 For $0<\ell<1$, we set
    $$ v^c_\ell(x)=T^1_\ell[u](x)=T^2_\ell[u](x)=\ell^{N-1}u(\ell x), \qq x\in \Gw^{\ell}=\frac{1}{\ell}\Gw. $$
By comparison principle, $v^c_\ell \leq P^{\Gw^\ell}(.,0)$ in $\Gw^\ell$  for every $\ell \in (0,1)$. Due to \rprop{halfspace-c}, up to a subsequence, $v^c_\ell$ converges to a function $v^c$ which is a solution of either \eqref{half1} if $H$ satisfies \eqref{multi}, or \eqref{half2} if $H$ satisfies \eqref{add} with $p=\frac{q}{2-q}$, or \eqref{half3} if $H$ satisfies \eqref{add} with $p>\frac{q}{2-q}$, or \eqref{half4} if $H$ satisfies \eqref{add} with $p<\frac{q}{2-q}$. Moreover, $v^c \leq P^{\BBR^N_+}(.,0)$ in $\BBR^N_+$. \medskip

If $H$ satisfies \eqref{add} with $p \neq \frac{q}{2-q}$ then since $m_{p,q}=p_c$, it follows from \cite{MV1} and \cite{NV} that $v^c=0$. \medskip

If $H$ satisfies \eqref{multi} or $H$ satisfies \eqref{add} with $p = \frac{q}{2-q}$ then set $\CV=\{ v: v$ is a solution of \eqref{A0} in $\BBR^N_+$, $v^c \leq v \leq P^{\BBR^N_+}(.,0) \}$ and put $\tl v:=\sup\CV$.

\noindent{\bf Assertion: } $\tl v$ is a solution of \eqref{halfspace} in $\BBR^N_+$. 

Indeed, let $\{Q_n\}$ be a sequence of $C^2$ bounded domain such that $\ovl Q_n \sbs Q_{n+1}$, $\cup_{n \in \BBN}Q_n=\BBR^N_+$ and $0<\dist(Q_n,\prt \BBR^N_+)<\frac{1}{n}$ for each $n \in \BBN$. Consider the problem
    \bel{Dn} \left\{ \BA{lll} -\Gd w +  H(x,w,\nabla w) & = 0 \qq &\text{in } Q_n \\
    \phantom{qerwerwhgh,,,,,,}
    w&= P^{\BBR^N_+}(.,0) &\text{on } \prt Q_n\EA \right. \ee
Since $v^c$ and $P^{\BBR^N_+}(.,0)$ are respectively subsolution and supersolution of \eqref{halfspace}, there exists a solution $w_n$ of the problem \eqref{halfspace} satisfying $0 \leq w_n \leq P^{\BBR^N_+}(.,0)$ in $Q_n$. Hence, by comparison principle $0 \leq w_{n+1} \leq w_n \leq P^{\BBR^N_+}(.,0)$ in $Q_n$ for each $n \in \BBN$. Therefore, $\tl w:=\lim_{n \to \infty}w_n \leq P^{\BBR^N_+}(.,0)$ in $\BBR^N_+$. Again, by \cite{Lib}, we obtain \eqref{lie} with $v_\ell$ replaced by $w_n$ and $\Gw^\ell$ replaced by $Q_n$. Thus $\tl w$ is a solution of \eqref{halfspace}. On one hand, by the definition of $\tl v$, $\tl w \leq \tl v$. On the other hand, $\tl v \leq w_n$ in $Q_n$ for every $n$, and consequently $\tl v \leq \tl w$ in $\BBR^N_+$. Thus $\tl v= \tl w$. \medskip

For every $\ell>0$, we set $w_\ell=T^1_\ell[\tl v]=T^2_\ell[\tl v]=\ell^{N-1}\tl v(\ell x)$ with $x \in \BBR^N_+$ then $w_\ell=\sup\CV$. Therefore $w_\ell=\tl v$ in $\BBR^N_+$ for every $\ell>0$. Hence $\tl v$ is self-similar, namely $\tl v$ can be written under the separable form
    $$ \tl v(r,\gs)=r^{N-1}\gw^*_i(\gs) \qq (r,\gs) \in (0,\infty) \ti S^{N-1}_+$$
where $\gw^*_i$ is the nonnegative solution of either \eqref{PH1} if $H$ satisfies \eqref{multi} or \eqref{PH2} if $H$ satisfies \eqref{add} with $p=\frac{q}{2-q}$. Since $H$ is critical, it follows from \rth{exis-uniq-SS} that $\gw^*_i \equiv 0$, hence $\tl v \equiv 0$. Thus $v^* \equiv 0$.

Hence
    \bel{t4} \lim_{n \to \infty}(\sup\{\abs{v^c_{\ell_n}(x)}+\abs{\nabla v^c_{\ell_n}(x)}: x \in \Gg_{R^{-1},R} \cap \Gw^{\ell_n}\})=0. \ee
Consequently,
    $$ \lim_{x \to 0}\abs{x}^{N-1}u(x)=0 \qq \text{and} \qq \lim_{x \to 0}\abs{x}^N\abs{\nabla u(x)}=0. $$
Therefore, $ \lim_{x \to 0}(\abs{x}^N\rho(x)^{-1}u(x))=0$, namely $u=o(P^\Gw(.,0))$. By comparison principle,$u \equiv 0$. \qed    \medskip

Finally, we deal with the case $q=2$.

\bth{q=2} Assume $q=2$.  If $u\in C(\overline\Gw\setminus\{0\})\cap C^2(\Gw)$ is a nonnegative solution of \eqref{A0}
 vanishing on $\prt\Gw\setminus\{0\}$ then $u\equiv 0$. \es
\Proof Put
$$v=\left\{ \BA{lll}1-e^{-\frac{1}{p+1}u^{p+1}} &\text{ if } H \text{ satisfies } (\ref{multi}), \\
1-e^{-u} &\text{ if } H \text{ satisfies } (\ref{add}) \text{ with } p=\frac{q}{2-q} \EA \right. $$
then $v \in C(\overline\Gw\setminus\{0\})\cap C^2(\Gw)$, $0 \leq v \leq 1$ and $v$ sattisfies
	\bel{v} -\Gd v \leq 0 \q \text{in } \Gw, \qq v=0 \q \text{on } \prt \Gw. \ee
Let $\eta_\gd$ be the solution of
	\bel{vd} -\Gd \eta_\gd = 0 \q \text{in } D_\gd, \qq \eta_\gd=v \q \text{on } \prt D_\gd \ee
then by the maximum principle  $v \leq \eta_\gd \leq 1$ in $D_\gd$. The sequence $\{\eta_\gd\}$ converges to an harmonic function $\eta^* \geq v$ as $\gd \to 0$. Since $0 \leq \eta^* \leq 1$ and $\eta^*=0$ on $\prt \Gw \sms \{0\}$, it follows that $\eta^* \equiv 0$. Hence $v \equiv 0$, so is $u$. \qed \medskip

\noindent{\bf Proof of Theorem F.} The proof follows immediately from \rcor{supercritical}, \rth{critical} and \rth{q=2}. \qed \bigskip

\noindent{\bf Acknowledgements} The authors were supported by the Israel Science Foundation founded by the Israel Academy of Sciences and Humanities, through grant 91/10.


{\small
\begin{thebibliography}{99}
%

%
\bibitem{AGMQ} S. Alarc\'on, J. Garc\'ia-Meli\'an and A. Quaas, {\em Keller-Osserman type conditions for some elliptic problems with gradient terms}, {\bf J. Differential Equations 252}, 886-914 (2012).
%
\bibitem{BaGi} C. Bandle and E. Giarrusso, {\em Boundary blow up for semilinear elliptic equations with nonlinear gradient terms}, {\bf Adv. Differential Equations 1},133–150 (1996).
%
%

%
%
\bibitem{BVi} M. F. Bidaut-V\'eron and L. Vivier, {\em An elliptic semilinear equation with source term involving boundary measures: the subcritical case}, {\bf Rev. Mat. Iberoamericana 16}, 477-513 (2000).
%
%
%
%
%
%
\bibitem{GT} D. Gilbarg and N. Trudinger, {\em Elliptic Partial Differential Equations of Second Order. Second edition}, Springer, Berlin (1983).
%
%
\bibitem{GV} A. Gmira and L. V\'eron, {\em Boundary singularities of solutions of some nonlinear elliptic equations}, {\bf Duke Math. J. 64}, 271-324 (1991).
%
%
\bibitem{Lg2} J. F. Le Gall, {\em A probabilistic approach to the trace at the boundary for solutions of a semilinear parabolic partial differential equation}, {\bf J. Appl. Math. Stochastic Anal. 9}, 399-414 (1996).
%
%
\bibitem{Lib} G. Liebermann, {\em Boundary regularity for solutions of degenerate elliptic equations}, {\bf Nonlinear Anal. 12}, 1203-1219 (1988).
%
\bibitem{LiLo} E. H. Lieb and M. Loss, {\em Analysis}, Grad. Stud. Math. 14, Amer. Math. Soc. (1997)
%
\bibitem{Li1} P. L. Lions, {\em Quelques remarques sur les probl\`eme elliptiques quasilineaires du second ordre}, {\bf J. Analyse Math. 45}, 234-254 (1985).
%
%
\bibitem{KaKr} J. L. Kazdan and R. J. Kramer, {\em Invariant criteria for existence of solutions to secondorder quasilinear elliptic equations}, {\bf Comm. Pure Appl. Math. 31}, 619-645 (1978).
%
\bibitem{La} O.A. Ladyzhenskaya and N. N.  Ural\'tseva,  {\em Linear and quasilinear elliptic equations}. Translated from the Russian by Scripta Technica, Inc. Translation editor: Leon Ehrenpreis Academic Press, New York-London (1968).
%
\bibitem{Ma} M. Marcus, {\em Complete classification of the positive solutions of $-\Gd u+u^q=0$}, {\bf Journal d'analyse Math.} (2012).
%
\bibitem{Mi} G. Mingione, {\em The Calder\'on-Zygmund theory for elliptic problems with measure data}, {\bf Ann. SNS Pica Cl. Sci. 6}, 195-261 (2007).
%
%
\bibitem{MV1} M. Marcus and L. V\'eron, {\em The boundary trace of positive solutions of semilinear elliptic equations: the subcritical case}, {\bf Arch. Rational Mech. Anal. 144}, 201-231 (1998).
%
\bibitem{MV2} M. Marcus and L. V\'eron, {\em The boundary trace of positive solutions of semilinear elliptic equations: the supercritical case}, {\bf J. Math. Pures Appl. (9) 77}, 481-524 (1998).
%
\bibitem{MV3} M. Marcus and L. V\'eron, {\em Removable singularities and boundary trace}, {\bf J. Math. Pures Appl. 80}, 879-900 (2001).
%
\bibitem{MV4} M. Marcus and L. V\'eron, {\em The boundary trace and generalized boundary value problem for semilinear elliptic equations with coercive absorption}, {\bf Comm. Pure Appl. Math. LVI}, 689-731 (2003).
%
\bibitem{MVbook} M. Marcus and L. V\'eron, {\em Nonlinear second order elliptic equations involving measures} (2013).
%
\bibitem{NV} P. T. Nguyen and L. V\'eron, {\em Boundary singularities of solutions to elliptic viscous Hamilton-Jacobi equations}, {\bf J. Funct. Anal. 263}, 1487-1538 (2012).
%
%
%
%
%
%
\bibitem{Tr1} N. Trudinger, {\em Local estimates for subsolutions and supersolutions of general second order elliptic quasilinear equations}, {\bf Invent. Math. 61}, 67-79 (1980).
%
%
\bibitem{V1} L. V\'eron, {\em Elliptic equations involving measures}, Stationary partial differential equations. Vol. I, 593--712, {\bf Handb. Differ. Equ.}, North-Holland, Amsterdam, 2004.
%
\bibitem{V2} L.V\'eron, {\em Singularities of solutions of second other Quasilinear Equations}, Pitman Research Notes in Math. Series 353, Adison Wesley, Longman 1996.
\end{thebibliography}}
\end{document}